\documentclass[12pt, a4paper]{amsart}

\usepackage[english]{babel}
\usepackage{amsfonts}
\usepackage{amsmath}
\usepackage{amssymb, latexsym}
\usepackage{graphicx}
\usepackage[usenames]{color}

\usepackage{fullpage}
\usepackage{amsmath}

\usepackage{mathrsfs}

\newtheorem{theorem}{Theorem}[section]
\newtheorem{proposition}[theorem]{Proposition}
\newtheorem{corollary}[theorem]{Corollary}
\newtheorem{lemma}[theorem]{Lemma}
\newtheorem{remark}[theorem]{Remark}
\newtheorem{definition}[theorem]{Definition}

\numberwithin{equation}{section} 

\def\calh{\mathcal{H}}
\def\calm{\mathcal{M}}

\def\calx{\mathcal{X}}

\def\bbe{{\mathbb{E}}}
\def\bbr{{\mathbb{R}}}
\def\bbp{{\mathbb{P}}}

\renewcommand{\Re}{{\rm Re} \,}
\renewcommand{\Im}{{\rm Im} \,}

\def\en t{{{\rm Z}\mkern-5.5mu{\rm Z}}}

\def\<{\left<}
\def\>{\right>}
\def\({\left(}
\def\){\right)}

\def\9{{\infty}}
\def\barr{\begin{array}}
\def\earr{\end{array}}

\def\wt{\widetilde}
\def\wh{\widehat}
\def\ol{\overline}

\def\vf{{\varphi}}
\def\lbb{{\lambda}}
\def\g{{\gamma}}
\def\a{{\alpha}}

\def\3{\subset }
\def\na{{\nabla}}

\def\ve{{\varepsilon}}
\def\p{{\partial}}

\begin{document}

\title{Scattering for stochastic nonlinear Schr\"odinger equations}

\author{Sebastian Herr}
\address[S. Herr and M. R\"ockner]{Fakult\"at f\"ur
Mathematik, Universit\"at Bielefeld, 33615 Bielefeld, Germany.}
\author{Michael R\"ockner}\author{Deng Zhang}
\address[D. Zhang]{Department of Mathematics,
Shanghai Jiao Tong University, 200240 Shanghai, China.}
\begin{abstract}
We study the scattering behavior of global solutions to stochastic
nonlinear Schr\"odinger equations with linear multiplicative noise.
In the case where the quadratic
variation of the noise is globally finite and the nonlinearity is defocusing,
we prove that the solutions scatter at infinity in the pseudo-conformal space
and in the energy space respectively,
including the energy-critical case.
Moreover, in the case where the noise is large, non-conservative and has infinite quadratic
variation, we show that the solutions scatter at infinity with high probability
for all energy-subcritical exponents.
\end{abstract}


\maketitle
\section{Introduction and main results} \label{Sec-Intro}

We are concerned with stochastic nonlinear Schr\"odinger equations with linear multiplicative noise and the long-time behaviour of their solutions. More precisely, we consider
\begin{equation} \label{equa-x}
\begin{split}
   idX&= \Delta Xdt   + \lbb |X|^{\a-1}X dt  -i\mu(t) Xdt + i\sum\limits_{k=1}^N  X  G_k(t)d\beta_k(t)
,  \\
   X(0)&=X_0.
 \end{split}
\end{equation}
Throughout this paper we assume that $\a>1$  and $d\geq 3$.
The choice $\lbb =-1$ (resp. $\lbb=1$) corresponds to the defocusing (resp. focusing) case \cite{SS99}.

The last term is taken in the sense of It\^o,
$\beta_k$ are real-valued
Brownian motions on a probability space $(\Omega, \mathscr{F},
\mathbb{P})$ with normal (in particular right-continuous) filtration
$(\mathscr{F}_t)_{t\geq 0}$,
$G_k(t)(x) := G_k(t,x)=g_k(t)\phi_k(x)$,
$g_k$ are real valued predictable processes,
$g_k \in L^2_{loc}(\bbr^+; \bbr)$ $\bbp$-a.s.,
and $\phi_k\in C^\9(\bbr^d; \mathbb{C})$.
$\mu$ is assumed to be of the form
$\mu(t,x) = \frac12 \sum_{k=1}^N |G_k(t,x)|^2$,
such that $t \mapsto |X(t)|_2^2$ is a continuous martingale.
In a quantum mechanical interpretation,
the $\beta_k$ are the output of the continuous measurement,
and
the probability measure
\begin{align*}
    \wh{\mathbb{P}}^T_{X_0} (d \omega)
    := (\mathbb{E}_{\mathbb{P}} [|X_0|_2^2])^{-1} |X(T,\omega)|_2^2\ \mathbb{P} (d \omega)
\end{align*}
is the physical probability law of the events occurring in $[0,T]$,
see
\cite{BG09} and references therein for more information.
In particular, when $\Re G_k =0$, $1\leq k\leq N$,
the mass $|X(t)|_2^2$ is pathwisely conserved,
thus the quantum system has a unitary evolution in a random environment.
This case also arises from molecular aggregates with thermal fluctuations,
we refer to \cite{BCIR94,BCIRG95} and references therein.

The global well-posedness of \eqref{equa-x} was first studied in \cite{BD99, BD03}
for more general linear multiplicative noise in the conservative case
for a restricted range of subcritical exponents.
For the full range of
mass- and energy-subcritical exponents,
the global well-posedness  is proved in
the recent papers  \cite{BRZ14, BRZ16}
in both the conservative and non-conservative cases,
based on the rescaling approach  and Strichartz estimates for lower order perturbations of the Laplacian.
Very recently, global well-posedness for
\eqref{equa-x} with quite general noise in the $L^2$ case is proved in \cite{H16}
by using the stochastic Strichartz estimates from \cite{BM14}.
In addition, we refer the reader to
\cite{BF12} for results on Schr\"odinger equations with potential perturbed by a temporal white noise,
to \cite{BM14,BHW17} for results on compact manifolds,
and to \cite{CG15, BD10, DT11} for Schr\"odinger equations with modulated dispersion.

In this paper, we are mainly concerned with the asymptotic behavior of solutions to
stochastic nonlinear Schr\"odinger equations.
More precisely, we focus on the scattering property of solutions,
which is of physical importance and, roughly speaking,
means that solutions behave asymptotically like those to linear Schr\"odinger equations.

There is an extensive literature on scattering in the deterministic case.
Let $\a(d)$ denote the Strauss exponent, see \eqref{ad} below.
In the defocusing case, scattering was proved in the pseudo-conformal space
(i.e., $\{u\in H^1: |\cdot| u(\cdot) \in L^2 \}$)
 if $\a\in[1+\a(d), 1+4/(d-2))$,
$d\geq 3$. For small initial data, scattering was also obtained
for $\a \in (1+ 4/d, 1+ 4/(d-2))$ in \cite{T85,CW92}.
Moreover, in the energy space $H^1$, the scattering
property of solutions was first proved
for the inter-critical case where $\a\in (1+4/d, 1+ 4/(d-2))$ in \cite{GV85}.
In the much more difficult energy-critical case $\a= 1+ 4/(d-2)$,
global well-posedness and scattering are established in \cite{CKSTT08} for $d=3$,
in \cite{RV07} for $d=4$, and in \cite{V07} for $d\geq 5$.
Moreover,  global well-posedness and scattering in $L^2$ are proved in \cite{D12} for the mass-critical exponent $\a= 1+4/d$, $d\geq 3$. In the focusing case, there exists a threshold for global well-posedness,
scattering and blow-up; we refer to \cite{FXC11, KM06, KV10, M15, M13} and references therein.

In the framework of stochastic mechanics,
developed by E. Nelson \cite{N85},
there are also several works devoted to potential scattering,
in terms of diffusions instead of wave functions.
See, e.g., \cite{C84,C85,S80}.

However, to the best of our knowledge,
there are few results on the scattering problem for
stochastic nonlinear Schr\"odinger equations \eqref{equa-x}.
One interesting question is that,
whether the scattering property
is preserved under thermal fluctuations?
Furthermore,
in the regime where the deterministic system fails to scatter,
will the input of some large noise have the effect to improve scattering with high probability?

One major challenge here lies in establishing global-in-time Strichartz estimates
for \eqref{equa-x},
which actually measure the dispersion and are closely
related to those of time dependent and long range perturbations of the Laplacian, see \eqref{equa-z} below.
It is known (see e.g. \cite{GVV06}) that Strichartz estimates may
fail for certain perturbations in the deterministic case.
Because of rapid fluctuations of Brownian motions at large time,
global-in-time Strichartz estimates for \eqref{equa-x}
in the case where $G_k(t,x)=\phi_k(x)$ (i.e. independent of time) are not expected to hold.
Moreover, another key difficulty arises from the failure of conservation of
several important quantities, such as the Hamiltonian and the pseudo-conformal energy.
Due to the complicated formulas of these quantities,
involving several stochastic integrations,
it is quite difficult to
derive a-priori estimates as in the deterministic case (e.g.\ the decay estimates of $|X(t)|_{L^{\a+1}}$ in the pseudo-conformal space).

In the present work, in the case where $\lbb=-1$ (defocusing) and the quadratic variation of the noise is globally bounded,
we prove scattering of global solutions to \eqref{equa-x}
in both the pseudo-conformal  and the energy space.
Precisely,
in the pseudo-conformal space,
we prove scattering for all $\a\in 1+\a(d), 1+4/(d-2))$, $d\geq 3$,
where $\a(d)$ is the Strauss exponent.
In the energy space,
we prove scattering for $\a\in[\max\{2,1+4/d\}, 1+4/(d-2))$, $3\leq d\leq 6$.
Moreover, in the energy-critical case,
assuming the global well-posedness of \eqref{equa-x}
we also obtain scattering in both spaces based on \cite{CKSTT08,RV07,V07}.
Thus, in dimensions $d=3,4$,
we obtain scattering in the energy space for all exponents $\a$ of the nonlinearity
ranging from the mass-critical exponent $1+4/d$
to the energy-critical exponent $1+4/(d-2)$.
We also mention that,
the mass-critical case is proved here based on the very recent work \cite{D12}.

Furthermore,
it is known from \cite[Theorem 7.5.2]{C03} that
the deterministic non-trivial solutions do not have any scattering states
in the case where $\a \in (1,1+2/d]$.
This motivates our further study on the regularization effect of the noise on scattering.
Actually, questions about the impact of noise on deterministic systems have attracted significant attention in the field of stochastic partial differential equations, see e.g. \cite{BRZ14.3, BD05, DFRV16, FGP10}.
Here, we prove that in the presence of a large
non-conservative spatially independent noise,
the stochastic solution to \eqref{equa-x} exists globally and scatters
at infinity with high probability
for all $\a\in(1,1+4/(d-2))$, both in the pseudo-conformal and the energy space.
This includes the range $(1,1+2/d]$
where scattering fails in the deterministic case.

We begin with recalling the definition of solutions to \eqref{equa-x}.
\begin{definition}\label{def-x}
Fix $T>0$. An $H^1$-solution to \eqref{equa-x} is an $H^1$-valued continuous
$(\mathscr{F}_t)$-adapted process $X=X(t)$, $t\in[0,T],$ such that $|X|^\a\in L^1([0,T], H^{-1})$
and it satisfies
\begin{align}\label{equa-x'}
   X(t) =&  X_0 -  \int^t_0 (i\Delta X(s) +\mu X(s) + \lambda i|X(s)|^{\alpha-1} X(s) )ds \nonumber \\
         & + \sum\limits_{k=1}^N\int ^t_0 X(s)G_k(s) d \beta_k(s),\ \forall t\in [0,T],
\end{align}
Here, the integral $\int ^t_0 X(s)G_k(s)d \beta_k(s)$ is taken in sense of It\^o,
and \eqref{equa-x'} is understood as an equation in $H^{-1}(\bbr^d)$.
\end{definition}

As in \cite{BRZ14,BRZ16},
we need the following assumption to assure global well-posedness of \eqref{equa-x}
in the energy-subcritical case discussed below.
\begin{itemize}
  \item[{\rm}(H0)]
  For each $1\leq k\leq N$, $0<T<\9$,
  $G_k(t,x)= g_k(t)\phi_k(x)$,
  $g_k$ are real valued predictable processes,
  $g_k \in L^\9(\Omega\times [0,T])$ and
  $\phi_k \in C^\9(\bbr^d,\mathbb{C})$ such that for any muti-index $\g$, $1\leq |\g|\leq 3$,
  \begin{align} \label{asymflat}
     \lim\limits_{|x|\to \9} |x|^2 |\partial_x^\g \phi_k(x)| = 0.
  \end{align}
\end{itemize}

\begin{theorem} \label{Thm-GWP}
Assume $(H0)$.
Let $1<\a<1+4/(d-2)_+$ if $\lbb=-1$, and $1<\a<1+4/d$  if $\lbb =1$.
For each $X_0\in H^1$ and $0<T<\infty$, there
exists a unique $H^1$-solution $X$ to (\ref{equa-x})
such that
\begin{align} \label{thmx1}
   X\in L^2(\Omega; C([0,T];H^1)) \cap L^{\a+1}(\Omega;
   C([0,T];L^{\a+1})),
\end{align}
and for any Strichartz pair $(\rho,\g)$ (see Section \ref{Sec-Stri}),
\begin{align}\label {thmx2}
  X \in L^\g(0,T;W^{1,\rho}),\ \ \mathbb{P}-a.s..
\end{align}
If $ X_0 \in \Sigma :=\{u\in H^1: |\cdot|u(\cdot) \in L^2\}$,
then $X\in L^2(\Omega; C([0,T];\Sigma))$, and for any Strichartz pair $(\rho,\g)$,
\begin{align}\label{thmx2*}
  \| |\cdot| X \|_{L^\g(0,T;L^{\rho})} < \9,\ \ \mathbb{P}-a.s..
\end{align}

Moreover, if in addition $\lbb =-1$ and
$g_k \in L^\9(\Omega; L^2(\bbr^+))$, $1\leq k\leq N$,
then for any $X_0 \in H^1$ and any $p\geq 1$,
\begin{align}  \label{glo-X-H1}
     \bbe  \sup\limits_{0<t<\9}  (|X(t)|^p_{H^1} + |X(t)|_{L^{\a+1}}^{p}) \leq C(p) <\9.
\end{align}

\end{theorem}
The proof is similar to that in \cite{BRZ14,BRZ16} and it is postponed to the Appendix.

In the defocusing energy-critical case, i.e., $\lbb=-1, \a=1+4/(d-2)_+$, the local well-posedness
for \eqref{equa-x} has been proved in \cite{BRZ16} for all $d\geq 1$, see also \cite{BD03} for $d\leq 5$.
However, the global well-posedness is much more difficult and remains still open.
One of the main difficulties is that several important quantities such as the Hamiltonian
are no longer conserved in the stochastic case.
In order to consider the scattering in the energy-critical case we a-priori assume  that
\begin{itemize}
  \item[{\rm}(H0')]
  In the case $\lbb=-1$, $\a=1+4/(d-2)$, for every $T>0$ and $X_0 \in H^1$
  there exists a unique $H^1$-solution $X$ to \eqref{equa-x}
  such that $X\in L^\g(0,T;W^{1,\rho})$, $\bbp$-a.s.,  for any Strichartz pair $(\rho,\g)$.
  In addition, if $X_0\in \Sigma$, then $\||\cdot| X\|_{L^\g(0,T;L^\rho)}<\9$, $\bbp$-a.s..
\end{itemize}

We first study the scattering property of global solutions in the pseudo-conformal space.
In this case,  the temporal functions $g_k$ in $(H0)$
are assumed to
satisfy suitable integrability
and to decay to zero with appropriate speed at infinity.
\begin{itemize}
\item[{\rm(H1)}] For each $1\leq k\leq N$,
\begin{align} \label{AF-3}
    \lim\limits_{|x|\to \9} |x|^3 |\partial_x^\g \phi_k(x)| = 0,\ \ 1\leq |\g|\leq 3,
\end{align}
${\rm esssup}_{\omega\in \Omega} \int_0^\9 (1+s^4) g_k^2(s) ds<\9$, and for $\bbp$-a.e. $\omega \in \Omega$,
\begin{align} \label{ILog}
     \lim\limits_{t\nearrow 1} (1-t)^{-3} \(\int_{\frac{t}{1-t}}^\9 g_k^2 (\omega, s)ds \ln \ln \({\int_{\frac{t}{1-t}}^\9 g_k^2 (\omega, s)ds}\)^{-1} \)^{\frac 12} =0.
\end{align}
\end{itemize}

Let $\a(d)$ denote the Strauss exponent, i.e.
\begin{equation} \label{ad}
\a(d) = \frac{2-d+\sqrt{d^2+12d+4}}{2d}.
\end{equation}

\begin{theorem}  \label{Thm-sca-X}
Let $X_0 \in \Sigma$, $\lbb=-1$ and
$\a\in (1+\a(d), 1+4/(d-2)]$, $d \geq 3$.
Assume $(H0)$, $(H1)$
and if $\a=1 + 4/(d-2)$ additionally $(H0')$.

$(i)$. $\bbp$-a.s.  there exists   $v_+\in \Sigma$, such that
\begin{align} \label{sca-X-delta}
    \lim\limits_{t\to \9} |e^{it \Delta}e^{-\vf_*(t)}X(t) - v_+|_{\Sigma} =0
\end{align}
with the rescaling function
\begin{align} \label{vf*}
   \vf_*(t)=-\sum_{k=1}^N \int_t^\9 G_k(s) d\beta_k(s) + \frac 12 \sum_{k=1}^N \int_t^\9 \(|G_k(s)|^2+G^2_k(s) \)ds.
\end{align}

$(ii).$
Let $V(t,s)$, $s,t\in [0,\9)$, be the evolution operators corresponding
to the random equation \eqref{equa-z*} below in the homogeneous case where $F \equiv 0$.
Then, $\bbp$-a.s. there exists  $X_+\in\Sigma$ such that
\begin{align} \label{sca-X-V}
    \lim\limits_{t\to \9} |V(0,t)e^{-\vf_*(t)} X(t) - X_+  |_{H^1} =0.
\end{align}
\end{theorem}

\begin{remark} \label{Rem-Integ-g-Sigma}
In Assumption $(H1)$, the $L^\9(\Omega)$-integrability of $\int_0^\9 (1+s^4) g_k^2(s)ds$
can be relaxed to the exponential integrability with respect to $\Omega$.
See Remark \ref{Rem-integ-g} below.
\end{remark}

The next result is concerned with the scattering in the energy space.
\begin{theorem} \label{Thm-Sca-H1}
Let $X_0 \in H^1$, $\lbb=-1$, $\a\in [\max\{2,1+4/d\}, 1+4/(d-2)]$, $d\geq3$.
Assume $(H0)$,  $g_k\in  L^2(0,\9)$, $\bbp$-a.s., $1\leq k\leq N$,
and if $\a=1+ 4/(d-2)$ additionally $(H0')$.

$(i)$. $\bbp$-a.s. there exists $v_+\in H^1$,  such that
\begin{align} \label{sca-H1}
    | e^{it\Delta} e^{-\vf_*(t)} X(t) - v_+ |_{H^1} \to 0,\ \ as\ t\to \9,
\end{align}
where $\vf_*$ is as in \eqref{vf*} above.

$(ii).$
Let $V(t,s)$, $s,t\in [0,\9)$, be the evolution operators as in Theorem \ref{Thm-sca-X}.
Then, $\bbp$-a.s., there exists  $X_+\in H^1$ such that
\begin{align} \label{sca-X-VH1}
    \lim\limits_{t\to \9}  |V(0,t)e^{-\vf_*(t)} X(t) - X_+  |_{H^1} =0.
\end{align}
\end{theorem}

\begin{remark}
One may remove the technical condition $\a\geq 2$ in Theorem \ref{Thm-Sca-H1} 
by using delicate arguments as in the $H^1$-critical stability result in \cite{TVZ07}. 
In order to keep the simplicity of exposition, 
we will not treat this technical problem in this paper.  
\end{remark}

Our next result is concerned with the regularization effect of noise on scattering
in the non-conservative case. We assume
\begin{enumerate}
  \item[(H2)] For each $1\leq k\leq N$, $\phi_k\equiv v_k$ are constants,
  $\inf_{t>0} g_k(t) \geq c_0>0$,
  and $\Re \phi_j \not =0$ for some $1\leq j\leq N$.
  Without loss of generality we may assume that $\Re \phi_1\not =0$.
\end{enumerate}

\begin{theorem} \label{Thm-sca-damped-X}
Let $X_0\in \Sigma$ (resp. $H^1$), $\lbb = \pm 1$, $\a \in (1,1+\frac{4}{d-2})$, $d\geq 3$.
Assume $(H0)$ and $(H2)$.
Let $g_k$ and $v_j$ being fixed, $1\leq k\leq N$, $2\leq j\leq N$,
and $A_{v_1}$ denote the event that
the solution $X$ to \eqref{equa-x} exists globally and scatters at infinity
in $\Sigma$ (resp. $H^1$),
namely, there exists a unique $u_+\in \Sigma$ (resp. $u_+\in H^1$) such that
\begin{align*}
    \lim\limits_{t\to \9} |e^{it\Delta} e^{-\vf(t)} X(t) - u_+|_{\Sigma} = 0 \ \
    (resp. \lim\limits_{t\to \9} |e^{it\Delta} e^{-\vf(t)} X(t) - u_+|_{H^1} = 0),
\end{align*}
where
\begin{align} \label{vf}
     \vf(t,x) := \sum\limits_{k=1}^N \int_0^t G_k(x,s) d\beta_k(s)
                 -\frac 12 \sum\limits_{k=1}^N \int_0^t  (|G_k(s,x)|^2 + G_k^2(s,x)) ds.
\end{align}
Then, we have
\begin{align} \label{sca-damped-X}
  \bbp (A_{v_1}) \to 1,\ as\ \Re v_1 \to \9.
\end{align}
\end{theorem}

\begin{remark}
We emphazise that
the scaling functions $\vf_*$ and $\vf$ depend on the strength of the noise (measured by the quadratic variation).
\end{remark}

\begin{remark}
It is known that for $\a\in (1,1+2/d]$,
scattering fails for solutions to deterministic nonlinear Schr\"odinger equations
(i.e., $G_k\equiv 0$).
Moreover, in the focusing case,
it is well-known that solutions may blow-up in the case where $\a\in [1+4/d,1+4/(d-2))$.
Hence, Theorem \ref{Thm-sca-damped-X} reveals a regularizing effect
of the noise on scattering.
See also \cite{BRZ14.3} for a regularizing effect of noise on blow-up
in the non-conservative case.
\end{remark}

\begin{remark}
It is interesting to consider the existence of wave operators and to raise the question
whether,
given any $v_+ \in \Sigma$ (resp. $H^1$),
there exists a unique solution to \eqref{equa-x}
such that the asymptotic behavior in Theorem \ref{Thm-sca-X} (resp. Theorem \ref{Thm-Sca-H1}) holds?
In the deterministic case,
the standard proof is to solve the equation backward in time,
that is, to first construct solutions on $[T, \9)$ for $T$ sufficiently large,
and then to extend the solution to all times.
However, the situation is quite different in the stochastic case.
Even if one can construct solutions path by path by using the deterministic strategy,
it is unclear whether the resulting solution is $\{\mathscr{F}_t\}$-adapted.
\end{remark}

Let us outline the proofs of Theorems \ref{Thm-sca-X}-\ref{Thm-sca-damped-X}. We rely on the rescaling approach
recently developed in \cite{BRZ14,BRZ16}
and on perturbative arguments.
By the transformation
\begin{align} \label{trans-z}
     z(t,x) := e^{-\vf(t,x)}X(t,x)
\end{align}
where $\vf$ is as \eqref{vf},
the original stochastic equation \eqref{equa-x} is reduced to the random equation below
\begin{align*}
   \partial_t z = -i e^{-\vf(t)} \Delta (e^{\vf(t)} z)  -\lbb i  e^{-\vf(t)} F (e^{\vf(t)} z),
\end{align*}
or equivalently,
\begin{align} \label{equa-z}
    \partial_tz =  A(t) z   -\lbb i e^{-\vf(t)} F (e^{\vf(t)} z).
\end{align}
Here, $F(u)=|u|^{\a-1}u$ for $u\in \mathbb{C}$,
$A(t):= -i(\Delta + b(t) \cdot \na + c(t))$, where
\begin{align}
  b(t) =& 2 \sum\limits_{k=1}^N \int_0^t \na G_k(s) d \beta_k(s) -2 \int_0^t \na \wh{\mu}(s)ds, \label{y-b}  \\
  c(t) =& \sum\limits_{j=1}^N \(\sum\limits_{k=1}^N \int_0^t \partial_j G_k(s) d \beta_k(s) - \int_0^t \partial_j \wh{\mu}(s)ds \)^2 \nonumber \\
        &  + \sum\limits_{k=1}^N\int_0^t \Delta G_k(s) d\beta_k(s) - \int_0^t \Delta \wh{\mu}(s)ds,  \label{y-c}
\end{align}
and
\begin{align*}
 \wh{\mu}(t,x) := \frac 12 \sum\limits_{k=1}^N (|G_k(t,x)|^2 + G_k(t,x)^2)
               = \sum\limits_{k=1}^N (\Re G_k(t,x)) G_k(t,x).
\end{align*}
(Note that, $\wh{\mu}=0$, if $\Re(G_k)=0$ for every $1\leq k \leq N$.
Otherwise, $\Re (\wh{\mu})>0$.)

Hence, the problem of scattering for \eqref{equa-x} is  reduced to that for the
random equation \eqref{equa-z}.

This point of view proved useful for a sharp pathwise analysis of stochastic solutions
and it reveals the structure of the initial stochastic equation as well.
Moreover, it is very  robust and applicable to several problems.
We refer to \cite{BR15,BRZ17}
for the applications of the rescaling approach,
combined with
the theory of maximal monotone operators,
to stochastic partial differential equations.
See also \cite{BRZ16.2} for optimal bilinear control problems
and \cite{Z17} for pathwise Strichartz
and local smoothing estimates for general stochastic dispersive equations.

The key observation here is that,
by \eqref{y-b} and \eqref{y-c},
the global bound  on the quadratic variation of the noise
implies that on the lower order coefficients,
which allows us to obtain the crucial
global-in-time Strichartz estimates for the time-dependent operator $A(t)$ in \eqref{equa-z},
see Theorem \ref{Thm-Stri} and Corollary \ref{Cor-Stri} below.
It should be mentioned that,
in the case where $G_k(t,x)= \mu_k e_k(x)$,
$\mu_k\in \mathbb{C}$, $e_k\in C^\9(\bbr^d)$,
as studied in \cite{BRZ14,BRZ16}, we have $\sup_{t\geq 0}(|b(t)|+|c(t)|)=\9$,
thus the global-in-time
Strichartz estimates  in this case are not expected to hold,
only local-in-time Strichartz estimates are available, see \cite{BRZ14,BRZ16,Z17}.

Heuristically, if $g_k\in L^2(\bbr^+)$,
the rescaling function $\vf$ converges almost surely at infinity, therefore
 the solution to
\eqref{equa-z} should behave asymptotically like that to the equation
\begin{align*}
   \partial_t z = -i e^{-\vf(\9)} \Delta (e^{\vf(\9)}z)  - \lbb i  e^{-\vf(\9)} F(e^{\vf(\9)}z).
\end{align*}
Thus,  $ z_*:=e^{\vf(\9)}z$ satisfies the equation
\begin{align*}
   \partial_t z_* = -i \Delta z_*  -\lbb i F(z_*),
\end{align*}
which is actually the deterministic nonlinear Schr\"odinger equation.
We shall mention that,
the papers \cite{C03,CW92,CKSTT08,RV07,V07} treat the equation
$\p_tu = i\Delta u - \lbb i |u|^{\a-1}u$, which
can be transformed to the equation above by reversing the time.
Thus, similar arguments apply also to this equation.
In particular, the solution $z_*$ scatters  at infinity
for appropriate $\alpha$.

Rigorously, we will use perturbative arguments to prove scattering properties of the
random solution $z$ to \eqref{equa-z}.
For this purpose,
we consider
\begin{align} \label{trans-z*}
z_* (t) :=e^{\vf(\9)}z(t)
         = e^{-(\vf(t) - \vf(\9))}X(t)
         = e^{-\vf_*(t)} X(t),
\end{align}
where $\vf$ and $\vf_*$ are as in \eqref{vf} and  \eqref{vf*}, respectively.
Then, it follows from \eqref{equa-z} that
\begin{align} \label{equa-z*}
     \partial_t z_* = A_*(t)z_* -\lbb i  e^{-\vf_*(t)} F(e^{\vf_*(t)}z_*),
\end{align}
where $A_*(t) = -i (\Delta + b_*(t) \cdot \na + c_*(t))$
with the coefficients
\begin{align}
     b_*(t) =& -2 \sum\limits_{k=1}^N \int_t^\9 \na G_k(s) d\beta_k(s) + 2  \int_t^\9 \na \wh{\mu}(s)ds, \label{b*}\\
     c_*(t) =& \sum\limits_{j=1}^N \( \sum\limits_{k=1}^N  \int_t^\9 \partial_j G_k(s) d\beta_k(s) - \int_t^\9 \partial_j\wh{\mu}(s)ds\)^2 \nonumber \\
               &-\sum\limits_{k=1}^N \int_t^\9 \Delta G_k(s) d\beta_k(s)  + \int_t^\9 \Delta \wh{\mu}(s)ds. \label{c*}
\end{align}
One important fact here is that the coefficients $b_*$, $c_*$
are asymptotically small at large time,
which suffices to yield global-in-time Strichartz estimates based on the work \cite{MMT08}, see Section \ref{Sec-Stri} below.

Next, in order to compare the solutions $z_*$ to \eqref{equa-z*} and $u$ to the
deterministic nonlinear Schr\"odinger equation
\begin{align} \label{equa-u}
   \p_t u = -i \Delta u   -\lbb i F(u)  ,
\end{align}
with  $u(T) = z_*(T)$ at large time $T$,
we set
\begin{align} \label{v}
   v:= z_*-u,
\end{align}
and obtain
\begin{align} \label{equa-v}
     & \partial_t v = A_*(t) v  -i (b_* \cdot \na + c_*) u
       -\lbb i (e^{-\vf_*(t)} F(e^{\vf_*(t)}(v+u)) - F(u)), \\
     & v(T) =0. \nonumber
\end{align}
At this step, the proof of Theorems \ref{Thm-sca-X} and \ref{Thm-Sca-H1}
is reduced to obtaining asymptotic estimates for the solution $v$ to \eqref{equa-v}.

For the scattering
in the pseudo-conformal space a key role
is played by an a-priori estimate in the scale-invariant space
$L^{\wt{q}}(0,\9; L^{\a+1})$ with $\wt{q} = \frac{2(\a^2 -1)}{4 -(\a-1)(d-2)}$,
which is proved in \cite{T85, C03} by the decay estimate of the $L^{\a+1}$-norm of solutions
from the pseudo-conformal conservation law
or by inhomogeneous Strichartz estimates for the non-admissible pair $(\a+1, \wt{q})$
\cite{CW92, M13, FXC11}.
However, due to the complicated nature of the evolution formula of the pseudo-conformal energy and the unavailability of
Strichartz and local smoothing estimates for the non-admissible pair $(\a+1, \wt{q})$, we
proceed differently. We shall compare solutions at the level of the pseudo-conformal transformations.
The method of pseudo-conformal transformations has been applied
successfully to prove scattering in the pseudo-conformal space, see e.g. \cite{C03,CW92,T85}.
It has advantage that
the scattering problem of the original solutions at infinity
is reduced to the Cauchy problem of their pseudo-conformal transformations at the singular time $1$,
which in turn can be analysed by Strichartz estimates
without relying on the decay estimates of the $L^{\a+1}$-norm of solutions.

Finally, the regularization effect of noise on scattering
in Theorem \ref{Thm-sca-damped-X}
can be  proved by the rescaling transformation as well.
The key point here is
that,
after the rescaling transformation,
one obtains an exponentially decaying term $e^{(\a-1)\Re \vf}$
in front of the nonlinearity,
which weakens the nonlinear effect
and allows random solutions  to scatter at infinity even in the case
$\a\in (1,1+2/d]$, $\lbb = \pm 1$,
where deterministic solutions fail to scatter.
From this perspective, the noise in the non-conservative case has a damping effect to the deterministic system.
We also refer the interested reader to
\cite{S15} for similar phenomena for deterministic damped fractional Schr\"odinger equations.

The remainder of this paper is organized as follows.
Sections  \ref{Sec-Finite-Varia} and \ref{Sec-H1} are mainly devoted to the proof of Theorems \ref{Thm-sca-X} and \ref{Thm-Sca-H1}.
Section \ref{Sec-Sca-Damp} is concerned with the proof of Theorem \ref{Thm-sca-damped-X}.
Finally, Section \ref{Sec-Stri} contains the global-in-time Strichartz and local smoothing estimates used in this paper.
The proofs of some auxiliary results are postponed to the Appendix.

\section{Scattering in the pseudo-conformal space} \label{Sec-Finite-Varia}

Recall that $\Sigma= \{v\in H^1: |\cdot|v(\cdot) \in L^2\}$ is the pseudo-conformal space.
Let $v$ be the solution to \eqref{equa-v} and
consider its pseudo-conformal transformation
\begin{align} \label{pct-wtv}
   \wt{v}(t,x) := (1-t)^{-\frac d2}\ v\(\frac{t}{1-t}, \frac{x}{1-t}\)\ e^{i\frac{|x|^2}{4(1-t)}},\ \ t\in [0,1),\ x\in \bbr^d.
\end{align}
By straightforward computations, we have
\begin{align} \label{equa-wtv}
   &\partial_t \wt{v}
   = \wt{A}_*(t) \wt{v} - i (\wt{b}_* \cdot \na + \wt{c}_*) \wt{u}
       - \lbb i h(t) (e^{-\wt{\vf}_*}F(e^{\wt{\vf}_*} (\wt{v}+\wt{u})) - F(\wt{u})),
\end{align}
with  $\wt{v}(\wt{T})=0$, $\wt{T} = T/(1+T) \in (0,1)$.
Here, $\wt{\vf}_*(t,x) = \vf_*(\frac {t}{1-t}, \frac {x}{1-t})$,
$h(t)= (1-t)^{\frac{d(\a-1)-4}{2}}$, and
$\wt{A}_*(t) = -i (\Delta + \wt{b}_*(t) \cdot \na + \wt{c}_*(t))$ with
\begin{align}
  \wt{b}_*(t,x) =&   (1-t)^{-1} b_* (\frac{t}{1-t}, \frac{x}{1-t}),  \label{wtb}\\
  \wt{c}_*(t,x) =&   (1-t)^{-2} \(-\frac i2   b_* (\frac{t}{1-t}, \frac{x}{1-t}) \cdot x
                    +  c_* (\frac{t}{1-t}, \frac{x}{1-t})\), \label{wtc}
\end{align}
where $b_*,c_*$ are as in \eqref{b*} and \eqref{c*} respectively.

Note that, under Assumptions $(H0)$ and $(H1)$,
global-in-time Strichartz estimates hold for the operator $\wt{A}_*$ on $[0,1)$.
See Section \ref{Sec-Stri} below for details.

Similarly, let $\wt{z}_*,\wt{u}$ denote the pseudo-conformal transformations of $z_*$ and $u$, respectively.
Then, we have
\begin{align} \label{equa-wtz}
  \partial_t \wt{z_*} = \wt{A}_*(t) \wt{z_*}  - \lbb i h(t) e^{-\wt{\vf}_*} F(e^{\wt{\vf}_*} \wt{z_*}),
\end{align}
and
\begin{align} \label{equa-wtu}
  \partial_t \wt{u} = -i \Delta \wt{u} -\lbb i h(t) F(\wt{u}).
\end{align}
with $\wt{u}(\wt{T}) = \wt{z}_*(\wt{T})$.
Note that, $\wt{u}$ depends on $\wt{T}$ and $\wt{z}_*(\wt{T})$.

For the solution $\wt{z_*}$ to \eqref{equa-wtz},
we first have the estimates on any bounded interval $[0,\wt{T}]$ below, $0<\wt{T}<1$.

\begin{lemma} \label{Lem-wtvT}
Assume the conditions in Theorem \ref{Thm-sca-X} to hold.
For each $\wt{v}(0) \in \Sigma$ and $0< \wt{T}<1$,
there exists a unique $H^1$-solution $\wt{z}_*$
to \eqref{equa-wtv} on $[0,\wt{T}]$, such that
for any Strichartz pair $(\rho, \g)$,
\begin{align} \label{glbdd-wtz-LpWq-Sigma}
   \|\wt{z}_*\|_{L^{\g}(0,\wt{T}; W^{1,\rho})}
   + \|\wt{z}_*\|_{LS{(0,\wt{T})}}
   + \||\cdot|\wt{z}_*\|_{L^{\g}(0,\wt{T}; L^{\rho})}
   + \|\wt{z}_*\|_{C([0,\wt{T}]; \Sigma)}
   <\9, \bbp-a.s..
\end{align}
\end{lemma}

{\it \bf  Proof.}
By direct computations, for any Strichartz pair $(\rho,\g)$ and $\wt{T}\in (0,1)$,
\begin{align} \label{wtz-z}
   &\|\wt{z}_*\|_{L^\g(0,\wt{T}; L^\rho)} =  \|z_*\|_{L^\g(0,\frac{\wt{T}}{1-\wt{T}}; L^\rho)}, \\
   &\|\na \wt{z}_*\|_{L^\g(0,\wt{T}; L^\rho)}
     \leq C (1+ \frac{\wt{T}}{1-\wt{T}}) \|z_*\|_{L^\g(0,\frac{\wt{T}}{1-\wt{T}}; W^{1,\rho})}
            + \||\cdot|z_*\|_{L^\g(0,\frac{\wt{T}}{1-\wt{T}}; L^\rho)}.
\end{align}
Thus, in view of \eqref{thmx2}, \eqref{thmx2*}, Assumption $(H0')$ and $\vf_* \in L^\9(0,\9;W^{1,\9})$,
we obtain that
$$ \|\wt{z}_*\|_{L^\g(0,\wt{T}; W^{1,\rho})}<\9,\ \  a.s..$$
In particular, $\|\wt{z}_*\|_{C([0,\wt{T}];H^1)} <\9$, a.s..

Moreover, we compute that for any $t\in [0,\wt{T}]$,
\begin{align*}
   |x_j \wt{z}_*(t)|_{2}^2
   =& \int (1-t)^{-d} |x_j z_*(\frac{t}{1-t}, \frac{x}{1-t})|^2 dx
   = \int (1-t)^{2} |y_j z_*(\frac{t}{1-t}, y)|^2 dy \\
   \leq& \|z_*\|^2_{C([0,\frac{\wt{T}}{1-\wt{T}}];\Sigma)} <\9,
\end{align*}
which implies that
$$\|\wt{z}_*\|_{C([0,\wt{T}];\Sigma)}<\9,\ \  a.s..$$

Now, for the estimate in the local smoothing space,
we take the Strichartz pair $(p,q)=(\frac{d(\a+1)}{d+\a-1}, \frac{4(\a+1)}{(d-2)(\a-1)})$.
Since $\a\in (1,1+4/(d-2)]$,
there exist $1<l<\9$, $1<\theta \leq \9$ such that
$1/q'=1/\theta+\a/q$ and $1/p'=1/l+1/p$.
Then,
taking into account $\|\wt{\vf}_*\|_{L^\9(0,1;W^{1,\9})}< \9$ a.s.,
applying   the Strichartz estimates to \eqref{equa-wtz} and using the H\"older inequality
we obtain
\begin{align*}
   \|\wt{z}_*\|_{LS(0,\wt{T})}
   \leq& C|\wt{z}_*(\wt{T})|_{2} + C \|h e^{-\wt{\vf}_*} F (e^{\wt{\vf}_*}\wt{z}_*)\|_{L^{q'}(0,\wt{T};L^{p'})} \\
   \leq& C|\wt{z}_*(\wt{T})|_{2}
         + C |h|_{L^\theta(0,\wt{T})} \|\wt{z}_*\|^{\a-1}_{L^q(0,\wt{T}; L^{(\a-1)l})} \|\wt{z}_*\|_{L^{q}(0,\wt{T};L^{p})}.
\end{align*}
Since
\begin{align} \label{h}
   h\in L^\theta(0,1), \ for\ \a\in [1+\a(d), 1+4/(d-2)] ,
\end{align}
taking into account the Sobolev imbedding $W^{1,p} \hookrightarrow L^{(\a-1)l}$ we obtain
\begin{align} \label{z*-lqp}
   \|\wt{z}_*\|_{LS(0,\wt{T})}
   \leq& C|\wt{z}_*(\wt{T})|_{2}
         + C |h|_{L^\theta(0,\wt{T})} \|\wt{z}_*\|^{\a}_{L^q(0,\wt{T}; W^{1,p})}<\9,\ \ a.s..
\end{align}

Concerning the weighted norm,
for each $1\leq j\leq d$, $x_j \wt{z}_*$ satisfies
\begin{align} \label{equa-xwtz}
   \p_t(x_j\wt{z}_*)
   = \wt{A}_*(t) (x_j\wt{z}_*) + i(2\p_j \wt{z}_* + \wt{b}_{*,j}\wt{z}_*) -\lbb i h(t) |e^{\wt{\vf}_*} \wt{z}_*|^{\a-1} (x_j \wt{z}_*).
\end{align}
Take a finite partition $\{t_k\}_{k=0}^M$ of $[0,\wt{T}]$
such that $\|z_*\|_{L^q(t_k,t_{k+1}; W^{1,p})} \leq \ve$,
and let  $\<x\> = \sqrt{1+|x|^2}$.
Arguing as above and using the fact that
\[\sup_{(t,x)\in [0,1]\times \bbr^d} \<x\>^2 |\wt{b}_*(t,x)| <\9,\]
(see Section \ref{Sec-Stri} below) we get
\begin{align} \label{lpq-xwtz}
   \|x_j \wt{z}_*\|_{L^q(t_k,t_{k+1}; L^p)}
   \leq& C\|\wt{z}_*\|_{C([0,\wt{T}]; \Sigma)}
       + C\|\p_j\wt{z}_* + \wt{b}_{*,j} \wt{z}_*\|_{L^1(t_k,t_{k+1}; L^2)} \nonumber \\
       &+ C \|he^{(\a-1)\Re \wt{\vf_*}} |\wt{z}_*|^{\a-1} (x_j\wt{z}_*)\|_{L^{q'}(t_k,t_{k+1}; L^{p'})} \nonumber \\
   \leq& C\|\wt{z}_*\|_{C([0,\wt{T}]; \Sigma)} + C \ve^{\a-1} \|x_j\wt{z}_*\|_{L^q(t_k,t_{k+1}; L^p)}.
\end{align}
Thus, taking $\ve$ small enough and then summing over $k$ we obtain
$\|x_j \wt{z}_*\|_{L^q(0,\wt{T}; L^p)}<\9 $, a.s..
Hence, the proof is complete.  \hfill $\square$\\

The following result,
involving the pseudo-conformal energy,
is crucial for the scattering behavior.
\begin{lemma}  \label{Lem-Ito-PC}
Let $\lbb=-1$,
$\a\in (1, 1+4/(d-2)]$, $d\geq 3$.
Assume $(H0)$, $(H0)'$ and  $(H1)$.
Let $X$ be the solution to \eqref{equa-x}
with $X(0)=X_0\in \Sigma$.
Define the pseudo-conformal energy
\begin{align} \label{PC-Energy}
  E(X(s)): = \int |(y-2i(1+s)\na) X(s,y)|^2 dy + \frac{8}{1+\a} (1+s)^2 |X(s)|_{L^{\a+1}}^{\a+1},\ \ s>0.
\end{align}
Then, we have $\bbp$-a.s.
\begin{align} \label{glo-Ito-E}
 E(X(s)) =& E(X_0) + \int_0^s a(r)dr
           + \frac{16}{\a+1} (1-\frac{d(\a-1)}{4}) \int_0^s (1+r) |X(r)|_{L^{\a+1}}^{\a+1} dr \nonumber \\
          &  + \sum\limits_{k=1}^N \int_0^s \sigma_k(r) d\beta_k(r),\ \ s>0,
\end{align}
where $a,\sigma_k$ are continuous $(\mathscr{F}_t)$-adapted processes, satisfying that
\begin{align} \label{bdd-asigma}
  \int_0^\9 |a(r)|dr + \sup\limits_{0\leq s<\9} \bigg|\sum\limits_{k=1}^N\int_0^s \sigma_k(r) d\beta_k(r)\bigg| <\9, \ \ a.s..
\end{align}
\end{lemma}

\begin{remark} \label{Rem-integ-g}
In Assumption $(H1)$,
one can relax
the $L^\9(\Omega)$-integrability of $\int_0^\9 (1+s^4)g_k^2(s)ds$
to
the weaker exponential integrability,
which is actually sufficient for
the almost surely global bound \eqref{bdd-asigma}
(see \eqref{esti-H-p} and \eqref{esti-asigma-g} below)
and for the pathwise estimates below as well.
\end{remark}

\subsection{Proof of Lemma \ref{Lem-Ito-PC}.}
Let $H(X), V(X), G(X)$ be the Hamiltonian, virial and momentum functions of $X$ as
in the proof of Theorem \ref{Thm-GWP} in the Appendix below, respectively.

Note that,
for the pseudo-conformal energy given by \eqref{PC-Energy},
\begin{align*}
   E(X(s))
   =8(1+s)^2 H(X(s)) - 4(1+s) G(X(s)) +  V(X(s)).
\end{align*}

Similarly to \cite[(4.11)]{BRZ14.3},
we have
\begin{align} \label{Ito-G}
   G(X(s))
   =& G(X_0)
     + 4 \int_0^s H(X(r))dr
     + \frac{4 \lbb}{\a+1} (1-\frac{d(\a-1)}{4}) \int_0^s|X(r)|_{L^{\a+1}}^{\a+1} dr \nonumber \\
    & + \int_0^s a_3(r) dr
     + \sum\limits_{k=1}^N \int_0^s \sigma_{3,k}(r) d\beta_{3,k}(r)
\end{align}
where
\begin{align*}
    a_3(r) =&  - \sum\limits_{k=1}^N \Im \int y\cdot \na G_k(r,y) |X(r,y)|^2 \ol{G}_k(r,y) dy, \nonumber \\
    \sigma_{3,k}(r) =&  d \int |X(r)|^2 \Im G_k(r,y) dy
     -2  \Im \int y \cdot \na X(r,y) \ol{X}(r,y) \ol{G}_k(r,y) dy.
\end{align*}

By the It\^o formulas of $V(X)$, $H(X)$, $G(X)$
in \eqref{Ito-V}, \eqref{hamil} and \eqref{Ito-G}, respectively,
\begin{align} \label{Ito-E}
     E(X(s))
   =& E(X_0) + \int_0^s a(r) dr
      - \frac{16 \lbb}{\a+1} (1-\frac{d(\a-1)}{4}) \int_0^s (1+r) |X(r)|_{L^{\a+1}}^{\a+1} dr \nonumber \\
    &  + \sum\limits_{k=1}^N \int_0^s \sigma_k(r) d\beta_k(r),
\end{align}
where
\begin{align} \label{Ito-E-a}
    a(r) = 8(1+r)^2 a_1(r) + 4(1+r) \sum\limits_{k=1}^N \Im \int y \na G_k(r,y) |X(r,y)|^2 \ol{G}_k(r,y) dy
\end{align}
with $a_1(r)$ as in \eqref{hamil} below and
\begin{align} \label{Ito-sigma}
    \sigma_k(r)
    =& 8(1+r)^2 \Re \<\na (G_k(r) X(r)), \na X(r)\>_2
     - 8\lbb (1+r)^2 \int \Re G_k(r,y) |X(r,y)|^{\a+1} dy   \nonumber \\
     &  - 4 (1+r) d  \int |X(r,y)|^2 \Im G_k(r,y) dy    \nonumber \\
     & +8 (1+r) \Im \int y \cdot \na X(r,y) \ol{X}(r,y) \ol{G}_k(r,y) dy
      +  2 \int |y|^2|X(r)|^2 \Re G_k(r,y) dy
\end{align}
with $\sigma_{1,k}, \sigma_{2,k}, \sigma_{3,k}$ as in \eqref{hamil}, \eqref{Ito-V} below and \eqref{Ito-G}, respectively.

Regarding \eqref{bdd-asigma},
we see that, by \eqref{Ito-E-a}, \eqref{Ito-sigma} and \eqref{glo-X-H1},
\begin{align} \label{esti-asigma-g}
    \bbe \int_0^\9 |a(s)| + |\sigma^2_k(s)| ds
    \leq& C \bbe \int_0^\9 (1+\calh(X(s)))^2 (1+s)^4g_k^2(s) ds \nonumber \\
    \leq& C \(\bbe \sup\limits_{0\leq s<\9} (1+\calh(X(s)))^4
            \bbe \(\int_0^\9 (1+s)^4 g_k^2(s) ds \)^2 \)^{\frac 12}
    <\9,
\end{align}
where $\calh(X) = \frac 12 |X|_{H^1}^2 + \frac{1}{\a+1} |X|_{L^{\a+1}}^{\a+1}$.
This implies \eqref{bdd-asigma} and so finishes the proof.
\hfill $\square$

As a consequence,
we have the crucial global bound for the solution $\wt{z}_*$ below.

\begin{corollary} \label{Cor-wtz-globH1}
Let $\lbb=-1$,
$\a\in (1+\a(d), 1+\frac{4}{d-2}]$ with $\a(d)$ as in \eqref{ad}, $d\geq 3$.
Assume $(H0)$, $(H0)'$ and  $(H1)$.
Then,
\begin{align} \label{glbdd-wtz-LpWq-Weight}
     \|\wt{z}_*\|_{L^{\gamma}(0,1;W^{1,\rho})}
     +\| |\cdot| \wt{z}_*\|_{L^{\gamma}(0,1;L^{\rho})}
      + \|\wt{z}_*\|_{LS{(0,1)}}
    + \|\partial_j \wt{z}_*\|_{LS{(0,1)}}
     < \9,\ \ a.s..
\end{align}
\end{corollary}

{\bf Proof of Corollary \ref{Cor-wtz-globH1}. }
We consider the cases $\a\in [1+4/d, 1+ 4/(d-2)]$
and $\a\in (1+\a(d), 1+4/d)$ separately below.

{\bf $(i)$ The case $\a\in [1+4/d, 1+ 4/(d-2)]$. }
Let $\wt{X}$ be the pseudo-conformal transformation of $X$.
Note that, if $t:= s/(1+s)$, $s\in [0,\9)$,
$$ E(X(s))
   = \wt{E}_1(\wt{X}(t))
  :=  4 |\na \wt{X}(t)|_2^2
      + \frac{8}{1+\a} (1-t)^{\frac d2 (\a-1) -2} |\wt{X}(t)|_{L^{\a+1}}^{\a+1},\ \ 0\leq t<1.$$
Then, using \eqref{Ito-E} and the fact that $|\wt{X}(t)|_{L^{\a+1}}^{\a+1} = (1-t)^{-d(\a-1)/2} |X(s)|_{L^{\a+1}}^{\a+1}$
we get that,
\begin{align} \label{Ito-wtE1}
    \wt{E}_1(\wt{X}(t))
   =& \wt{E}_1(\wt{X}_0) + \int_0^sa(r)dr + \sum\limits_{k=1}^N \int_0^s \sigma_{k}(r)d\beta_k(r) \nonumber \\
    & + \frac{16}{\a+1} (1- \frac{d}{4}(\a-1))
              \int_0^t (1-r)^{\frac d2 (\a-1) -3} |\wt{X}(r)|_{L^{\a+1}}^{\a+1} dr     \\
   \leq& C|X_0|_{\Sigma}^2
        + \sup\limits_{0\leq s<\9}\( \ \bigg|\int_0^sa(r)dr\bigg|
        + \sum\limits_{k=1}^N \bigg| \int_0^s \sigma_{k}(r)d\beta_k(r) \bigg| \)<\9,\ \ a.s.,  \nonumber
\end{align}
where the last step is due to \eqref{bdd-asigma}.
Since
$ |\na \wt{X}(t)|_2^2 \leq \frac 14   \wt{E}_1(\wt{X}(t))$
and by \eqref{glo-X-H1},
$\sup_{t\in[0,1)}|\wt{X}(t)|_2^2 = \sup_{s\in[0,\9)} |X(s)|_2^2 <\9$, a.s.,
we obtain that
$ \sup_{t\in[0,1)} |\wt{X}(t)|_{H^1} <\9$, a.s..

Thus,
taking into account that $\wt{z}_* = e^{\wt{\vf}_*} \wt{X}$
and $\wt{\vf}_* \in L^\9(0,1; W^{1,\9})$,
we get
\begin{align} \label{bdd-wtz-H1-case1}
    \sup\limits_{t\in[0,1)} |\wt{z}_*(t)|_{H^1} <\9, \ \ a.s..
\end{align}

We claim that the estimate \eqref{bdd-wtz-H1-case1} is sufficient to yield  that
for $\wt{T}$ close to $1$,
\begin{align} \label{esti-wtz-LpWq-Weight-case1}
   \|\wt{z}_*\|_{L^{\gamma}(\wt{T},1;W^{1,\rho})}
     +\| |\cdot| \wt{z}_*\|_{L^{\gamma}(\wt{T},1;L^{\rho})}
   \leq  C(|\wt{z}_*(\wt{T})|_{\Sigma}  +  \|\wt{z}_*\|_{C([0,1);H^1)})
   <\9.
\end{align}
To this end,
we choose the Strichartz pair $(p,q)$, $1<\theta\leq \9$ and $h$ as in the proof of Lemma \ref{Lem-wtvT}.
Similarly to \eqref{z*-lqp}, we have that for any $t\in (\wt{T}, 1)$,
\begin{align*}
      \|\wt{z}_*\|_{L^{q}(\wt{T},t;W^{1,p})}
      \leq& C |\wt{z}_*(\wt{T})|_{H^1} + C |h|_{L^\theta(\wt{T}, t)}
            \|\wt{z}_*\|^\a_{L^{q}(\wt{T},t;W^{1,p})}  \\
      \leq& C \|\wt{z}_*\|_{C([0,1);H^1)}
           + C \ve(\wt{T})
            \|\wt{z}_*\|^\a_{L^{q}(\wt{T},t;W^{1,p})},
\end{align*}
where in the last step we also used
that  $\ve(\wt{T}) := |h|_{L^\theta(\wt{T}, 1)} \to 0$ as $\wt{T} \to 1^{-}$
(see also \eqref{h-0-Sigma} below),
due to \eqref{h-0-Sigma}.
Taking $\wt{T}$ close to $1$,
using \cite[Lemma A.1]{BRZ16.2}
and then letting $t\to 1^{-}$
we obtain
\begin{align} \label{esti-wtz-LpWq-case1}
\|\wt{z}_*\|_{L^{q}(\wt{T},1;W^{1,p})}
    \leq C \|\wt{z}_*\|_{C([0,1);H^1)} ,\ \ a.s..
\end{align}
Moreover, similarly to \eqref{lpq-xwtz}, for each $1\leq j\leq d$,
\begin{align*}
   \|x_j \wt{z}_*\|_{L^q(\wt{T}, t; L^p)}
   \leq C|\wt{z}_*(\wt{T})|_{\Sigma}
       +C  \|\wt{z}_*\|_{C([0,1);H^1)}
       + C \ve(\wt{T}) \|x_j\wt{z}_*\|_{L^q(\wt{T}, t; L^p)}.
\end{align*}
Then, taking $\wt{T}$ close to $1$
such that $C\ve(\wt{T}) \leq 1/2$
and letting $t\to 1^-$ we obtain
\begin{align} \label{esti-wtz-weight-case1}
   \|x_j \wt{z}_*\|_{L^q(\wt{T}, 1; L^p)}
   \leq C(|\wt{z}_*(\wt{T})|_{\Sigma}  +  \|\wt{z}_*\|_{C([0,1);H^1)}) .
\end{align}
Thus, combining \eqref{esti-wtz-LpWq-case1} and \eqref{esti-wtz-weight-case1} together
and using Strichartz estimates
we prove \eqref{esti-wtz-LpWq-Weight-case1}, as claimed.
Now, taking into account Lemma \ref{Lem-wtvT}
and using Strichartz estimates to control the $LS(0,1)$-norms
we obtain \eqref{glbdd-wtz-LpWq-Weight} in the case
where $\a\in [1+4/d, 1+4/(d-2)]$.

{\bf The case $\a\in (1+\a(d), 1+4/d)$.} In this case,
let
$$ \wt{E}_2(\wt{X}(t)) := (1-t)^{2-\frac d2 (\a-1)} \wt{E}_1(\wt{X}(t))
   =4(1-t)^{2-\frac d2 (\a-1)} |\na\wt{X}(t)|_2^2 + \frac{8}{1+\a} |\wt{X}(t)|_{L^{\a+1}}^{\a+1}. $$
Note that, by \eqref{Ito-wtE1},
\begin{align}  \label{Ito-wtE2}
  \wt{E}_2(X(t))
  =& \wt{E}_2(X_0) + \int_0^s (1+r)^{\frac d2 (\a-1)-2} a(r) dr
    + \sum\limits_{k=1}^N \int_0^s (1+r)^{\frac d2 (\a-1)-2} \sigma_k(r) d\beta_k(r) \nonumber  \\
  & -8  (1- \frac{d}{4}(\a-1)) \int_0^t (1-r)^{1 - \frac d2 (\a-1)} |\na \wt{X}(r)|_2^2 dr  \\
  \leq& C|X_0|_\Sigma^2
        + \sup\limits_{0\leq s<\9} \(\bigg|\int_0^sa(r)dr\bigg|
        + \sum\limits_{k=1}^N \bigg| \int_0^s \sigma_{k}(r)d\beta_k(r) \bigg|\)  <\9,\ \ a.s.. \nonumber
\end{align}
This yields  that
${\rm sup}_{0\leq t<1} |\wt{X}(t)|_{L^{\a+1}}^{\a+1} <\9$,
and so
\begin{align} \label{bdd-wtz-Lp-case2}
  \sup_{0\leq t<1} |\wt{z}_*(t)|_{L^{\a+1}}^{\a+1} <\9, \ \ a.s..
\end{align}

As in the previous case, we claim that
the estimate \eqref{esti-wtz-LpWq-Weight-case1}
also holds in the case where $\a\in(1+\a(d), 1+4/d)$.

To this end, we choose the Strichartz pair $(p,q)= (\a+1, \frac{4(\a+1)}{d(\a-1)})$
and set $\wt{q} = \frac{2(\a^2-1)}{4-(\a-1)(d-2)}$.
Note that, since $\a>1+\a(d)$,
$h^{\frac {1}{\a-1}}\in L^{\wt{q}}(0,1)$,
and so
$ \ve'(\wt{T}) := |h^{\frac{1}{\a-1}} |^{\a-1}_{L^{\wt{q}}(\wt{T},1; L^p)} \to 0$ as
$\wt{T} \to 1^{-}$.
Applying Strichartz estimates to \eqref{equa-z*}
and using H\"older's inequality and \eqref{bdd-wtz-Lp-case2} we have
for any $t\in (\wt{T}, 1)$
\begin{align*}
    \|\wt{z}_*\|_{L^{q}(\wt{T},t;W^{1,p})}
    \leq& C|\wt{z}_*(\wt{T})|_{H^1}
       + C \|h^{\frac{1}{\a-1}} z_*\|^{\a-1}_{L^{\wt{q}}(\wt{T},t; L^p)}
         \|\wt{z}_*\|_{L^q(\wt{T},t; W^{1,p})} \\
   \leq& C|\wt{z}_*(\wt{T})|_{H^1}
        + C \ve'(\wt{T})
         \|\wt{z}_*\|_{L^q(\wt{T},t; W^{1,p})}.
\end{align*}
Then, taking $\wt{T}$ close to $1$ such that $\ve'(\wt{T}) = 1/(2C)$
and then letting $t\to 1^{-}$,
we obtain
\begin{align} \label{esti-wtz-LpWq-case2}
   \|\wt{z}_*\|_{L^{q}(\wt{T},1;W^{1,p})}
  \leq  C|\wt{z}_*(\wt{T})|_{H^1}.
\end{align}
This, via Strichartz estimates, yields that
$\|\wt{z}_*\|_{C([\wt{T},1);H^1)} \leq C|\wt{z}_*(\wt{T})|_{H^1}<\9$, a.s..

Similarly,
for the estimate in the weighted space,
we get that for each $1\leq j\leq d$,
similarly to \eqref{lpq-xwtz},
for any $t\in(\wt{T}, 1)$,
\begin{align*}
   \|x_j \wt{z}_*\|_{L^q(\wt{T}, t; L^p)}
   \leq& C |\wt{z}_*(\wt{T})|_{\Sigma}
         + C \|\wt{z}_*\|_{C([\wt{T}, 1); H^1)}
         +  C \|h^{\frac{1}{\a-1}} z_*\|^{\a-1}_{L^{\wt{q}}(\wt{T},t; L^p)}  \|x_j\wt{z}_*\|_{L^q(\wt{T}, t; L^p)} \\
   \leq& C |\wt{z}_*(\wt{T})|_{\Sigma}
         + C \|\wt{z}_*\|_{C([\wt{T}, 1); H^1)}
         + C \ve'(T) \|x_j\wt{z}_*\|_{L^q(\wt{T}, t; L^p)}.
\end{align*}
Thus, similar arguments as above yield
that
\begin{align} \label{esti-wtz-Weight-case2}
   \|x_j \wt{z}_*\|_{L^q(\wt{T}, 1; L^p)}
   \leq& C\|\wt{z}_*(\wt{T})\|_{\Sigma}
         + C \|\wt{z}_*\|_{C([\wt{T}, 1); H^1)}
   <\9, \ a.s..
\end{align}

Now, we can use \eqref{esti-wtz-LpWq-case2},
\eqref{esti-wtz-Weight-case2}
and Strichartz estimates to
obtain \eqref{esti-wtz-LpWq-Weight-case1}
in the case where $\a\in (1+\a(d), 1+4/d)$, as claimed.
Then, similar arguments as those below \eqref{esti-wtz-weight-case1}
yields \eqref{glbdd-wtz-LpWq-Weight}.

Therefore, the proof is complete.
\hfill $\square$ \\

It should be mentioned that,
on the basis of Lemma \ref{Lem-Ito-PC},
one can prove the scattering of $X$
by using  the Galilean operator and decay estimates as in \cite{C03, TVZ07}.
Moreover, exploring the global bound \eqref{glbdd-wtz-LpWq-Weight}
and the equivalence  of asymptotics of $\wt{z}_*$ at time $1$
and $z_*$ at infinity (see e.g. \cite[Proposition $7.5.1$]{C03}),
one can also obtain the scattering  \eqref{sca-X-delta}.
Precisely,
almost surely there exists $v_+ \in \Sigma$ such that
\begin{align*}
   \lim\limits_{t\to \9} |e^{it \Delta}z_*(t) - v_+|_{\Sigma}
  = \lim\limits_{t\to \9} |e^{it \Delta}e^{-\vf_*(t)}X(t) - v_+|_{\Sigma} =0.
\end{align*}

We remark that, though Proposition 7.5.1 in \cite{C03} treats the equation
$\p_t u = i \Delta u - \lbb i |u|^{\a-1} u$, similar arguments apply also to
the solutions $z_*$ to \eqref{equa-z*} and
$\wt{z}_*$ to \eqref{equa-wtz} considered here.
In fact, define the dilation $D_\beta$ by $D_\beta z_*(x) = \beta^{\frac d2} z_*(\beta x)$,  $\beta>0$,
the multiplication $M_\sigma$ by $M_\sigma z_*(x)= e^{i\frac{\sigma |x|^2}{4}} z_* (x)$, $\sigma\in \bbr$,
and let
$\mathcal{T}(t):= e^{it\Delta}$, $t\in \bbr$.
We have
$\wt{z}_*(t) = M_{\frac{1}{1-t}} D(\frac{1}{1-t}) z_*(\frac{t}{1-t})$, $t\in [0,1)$.
Since
$\mathcal{T}(t) D_\beta = D_\beta \mathcal{T}(\beta^2t)$,
$\mathcal{T}(t) M_\sigma = M_{\frac{\sigma}{1+\sigma t}} D_{\frac{1}{1+\sigma t}} \mathcal{T}(\frac{t}{1+\sigma t})$,
we obtain $\mathcal{T}(t) \wt{z}_*(t) = M_{1} \mathcal{T}(\frac{t}{1-t}) z_*(\frac{t}{1-t})$.
It follows that
$\mathcal{T}(s) z_*(s) = M_{-1} \mathcal{T}(\frac{s}{1+s}) \wt{z}_*(\frac{s}{1+s})$, $s\in [0,\9)$,
which implies the equivalence of asymptotics between $\wt{z}_*$ and $e^{it\Delta}z_*$.  \\

However, in the $H^1$ case, one can not obtain the scattering behavior of $z_*$
directly from the uniform bound similar to \eqref{glbdd-wtz-LpWq-Weight}.
Thus, we present a different proof for scattering to illustrate the idea of comparison,
which applies also to the $H^1$ case.
Moreover, it also gives the asymptotical estimates of the solutions $\wt{z}_*$ to \eqref{equa-wtz}
and $\wt{u}$  to \eqref{equa-wtu} (see \eqref{asym-wtv-Sigma} below)
and so justifies the intuition mentioned in Section \ref{Sec-Intro}.

Proposition \ref{Prop-wtu} below summarizes uniform estimates (independent of $\wt{T}$)
for $\wt{u}$ used in this section.
\begin{proposition} \label{Prop-wtu}
Let $\lbb=-1$,
$\a\in (1+\a(d), 1+\frac{4}{d-2}]$ with $\a(d)$ as in \eqref{ad}, $d\geq 3$.
Then,
for each $\wt{z}_*(\wt{T})\in \Sigma$,
there exists a unique  $H^1$-solution $\wt{u}$
(depending on $\wt{T}$) to \eqref{equa-wtu} on $[0,1]$
with $\wt{u}(\wt{T}) = \wt{z}_*(\wt{T})$,
such that $\bbp$-a.s.
$\wt{u}\in C([0,1];\Sigma)$,
and  for  any Strichartz pair $(\gamma,\rho)$,
\begin{align} \label{glbdd-wtu-LpWq-Weight}
    \|\wt{u}\|_{L^{\gamma}(\wt{T},1;W^{1,\rho})}
     +\| |\cdot| \wt{u}\|_{L^{\gamma}(\wt{T},1;L^{\rho})}
      + \|\wt{u}\|_{LS{(\wt{T},1)}}
    + \|\partial_j \wt{u}\|_{LS{(\wt{T},1)}}
    \leq C  < \9,
\end{align}
where $1\leq j\leq d$,
$LS(0,1)$ is the local smoothing space defined in Section \ref{Sec-Stri},
and $C$ is independent of $\wt{T}$.
\end{proposition}
The proof
is similar to that of Corollary \ref{Cor-wtz-globH1},
based on the pseudo-conformal energy
and the global bound \eqref{glbdd-wtz-LpWq-Weight} of $\wt{z}_*(\wt{T})$.
For simplicity, it is postponed to the Appendix.

The next lemma contains the crucial global estimates
and asymptotics  of the solution $\wt{v}$ to equation \eqref{equa-wtv}.
\begin{lemma} \label{Lem-h1-wtv}
Assume the conditions in Theorem \ref{Thm-sca-X} to hold.
Let $\wt{v}$ be the solution to \eqref{equa-wtv}
with $\wt{v}(\wt{T}) = 0$.
Then,
for any Strichartz pair $(\rho,\gamma)$,
\begin{align} \label{glbdd-wtv-Sigma}
    \|\wt{v}\|_{L^{\g}(\wt{T},1;W^{1,\rho})}
    + \|\wt{v}\|_{LS{(\wt{T},1)}}
    + \||\cdot| \wt{v}\|_{L^{\g}(\wt{T},1;L^{\rho})}
    + \| \wt{v}\|_{C([\wt{T},1);\Sigma)}
    \leq C < \9,\ \ a.s.,
\end{align}
where $C$ is independent of $\wt{T}$.
Moreover,  $\bbp$-a.s., as $\wt{T} \to 1^-$,
\begin{align} \label{asym-wtv-Sigma}
     \|\wt{v}\|_{L^{\g}(\wt{T},1; W^{1,\rho})}
   + \|\wt{v}\|_{LS{(\wt{T},1)}}
   + \||\cdot|\wt{v}\|_{L^{\g}(\wt{T},1; L^{\rho})}
   + \|\wt{v}\|_{C([\wt{T},1); \Sigma)}
   \to 0.
\end{align}
\end{lemma}

{\it \bf Proof.}
The uniform bound \eqref{glbdd-wtv-Sigma} follows
from \eqref{glbdd-wtz-LpWq-Weight}
and \eqref{glbdd-wtu-LpWq-Weight},
since $\wt{v} = \wt{z}_* -\wt{u}$.

Regarding \eqref{asym-wtv-Sigma},
the proof is based on perturbative arguments.
Let $p,q,\theta,l$ be as in the proof of Lemma \ref{Lem-wtvT}.
First, we consider the estimates in the spaces $L^\g(0,1;W^{1,\rho})$
and $LS(0,1)$,
where $(\rho, \g)$ is any Strichartz pair.
Similarly to \eqref{z*-lqp}, since $\wt{v}(\wt{T})=0$, by \eqref{equa-wtv} we have $\bbp$-a.s.
for any $t\in (\wt{T}, 1)$,
\begin{align} \label{esti-wtvf}
       & \|\wt{v}\|_{L^\g(\wt{T},t;L^\rho)} + \|\wt{v}\|_{C([\wt{T},t];L^2)}
         +  \| \wt{v}\|_{LS{(\wt{T},t)}}  \nonumber \\
    \leq&  C \|(\wt{b}_* \cdot \na+\wt{c}_*) \wt{u}\|_{LS'{(\wt{T},t)}}
         + C|h|_{L^\theta(\wt{T},t)}(\|\wt{u}\|^\a_{L^q(\wt{T},t;W^{1,p})}
         + \|\wt{v}\|^\a_{L^q(\wt{T},t;W^{1,p})})
\end{align}
where  $C$ is independent of $\wt{T}$ and $t$,
$LS'(\wt{T},t)$ is the dual space of $LS(\wt{T},t)$.

Taking into account \eqref{h} and that $|h|_{L^\9(\wt{T},t)} \leq (1-\wt{T})^{4/(d-2)}$
if $\a=1+4/(d-2)$,
we have that
\begin{align} \label{h-0-Sigma}
\ve_1(\wt{T}) := |h|_{L^\theta(\wt{T},1)}  \to 0\ \ as\ \wt{T} \to 1^{-}.
\end{align}

Moreover,
by Assumption $(H1)$,
for any $0\leq |\beta| \leq 2$, $0\leq |\g|\leq 1$,
$\p_x^\beta \wt{b}_*$ and $\p_x^\g \wt{c}_*$ satisfy
\begin{align}
   &\lim\limits_{t\to 1} \sup\limits_{\bbr^d}
    \<x\>^2 (|\partial_x^\beta \wt{b}_*(t,x)| + \partial_x^\g \wt{c}_*(t,x)|) =0, \ \ a.s., \label{esti-wtvf.1}\\
   &\lim_{|x|\to \9}  \sup\limits_{t\in[0,1)}
    \<x\>^2 (|\partial_x^\beta \wt{b}_*(t,x)| + \partial_x^\g \wt{c}_*(t,x)|) =0. \ \ a.s.. \label{esti-wtvf.2}
\end{align}
(See Section \ref{Sec-Stri} below.)
Then, using \eqref{X'-X} below (see also Remark \ref{Rem-Stri})
and the uniform bound \eqref{glbdd-wtu-LpWq-Weight}
we get $\bbp$-a.s.
\begin{align} \label{esti-b}
   \|(\wt{b}_* \cdot \na + \wt{c}_*) \wt{u}\|_{LS'{(\wt{T}, t)}}
   \leq \ve_2(\wt{T})  \|\wt{u}\|_{LS{(\wt{T}, t)}},
\end{align}
where $\ve_2(\wt{T}) \to 0$ as $\wt{T} \to 1^-$.

Thus, setting $\ve(\wt{T}) := \ve_1(\wt{T}) \vee \ve_2(\wt{T})$
and using \eqref{esti-wtvf}, \eqref{h-0-Sigma} and \eqref{esti-b} we get
\begin{align} \label{ho-wtv*}
    &\|\wt{v}\|_{L^\g(\wt{T},t;L^{\rho})} + \|\wt{v}\|_{C([\wt{T},t];L^2)} + \|\wt{v}\|_{LS{(\wt{T}, t)}} \nonumber \\
    \leq& C \ve(\wt{T})
          ( \|\wt{u}\|_{LS{(\wt{T},t)}} +  \|\wt{u}\|^\a_{L^q(\wt{T},t;W^{1,p})}
            +  \|\wt{v}\|^\a_{L^q(\wt{T},t;W^{1,p})} ) \nonumber \\
    \leq& C \ve(\wt{T}),
\end{align}
where $C$ is independent of $\wt{T}$ and $t$,
due to the uniform bounds \eqref{glbdd-wtu-LpWq-Weight}
and \eqref{glbdd-wtv-Sigma}.

Furthermore, for every $1\leq j\leq d$,  by \eqref{equa-wtv},
$\p_j\wt{v}$ satisfies the equation
\begin{align*}
    \partial_t(\partial_j \wt{v})
    =&\wt{A}_*(t) (\partial_j \wt{v})  + G( \wt{v}, \partial_j \wt{u}, \wt{u})
\end{align*}
where $\p_j\wt{v}(\wt{T}) = 0$, and
\begin{align*}
     G(\wt{v}, \partial_j \wt{u}, \wt{u})
     =& - i (\wt{b}_* \na + \wt{c}_* ) (\partial_j \wt{u})
     - i (\partial_j \wt{b}_* \cdot \na +\partial_j \wt{c}_* ) (\wt{v} + \wt{u})   \\
      & -\lbb i h(t) (\partial_j(e^{-\wt{\vf}_*}F(e^{\wt{\vf}_*}(\wt{v}+\wt{u}))) - \partial_j F(\wt{u})).
\end{align*}
Then, similar to \eqref{esti-wtvf}, we have $\bbp$-a.s.
\begin{align*}
    & \|\partial_j \wt{v}\|_{L^\g(\wt{T},t;L^\rho)} + \| \partial_j \wt{v}\|_{C([\wt{T},t];L^2)} \nonumber  \\
    \leq& C(\|(\wt{b}_* \cdot \na + \wt{c}_*) \partial_j \wt{u}
            + (\partial_j \wt{b}_* \cdot \na + \partial_j\wt{c}_*)(\wt{v} + \wt{u})\|_{LS'{(\wt{T},t)}})  \\
        & + C |h|_{L^\theta(\wt{T}, t)}  (\|\wt{v}\|^\a_{L^q(\wt{T},t;W^{1,p})} + \|\wt{u}\|^\a_{L^q(\wt{T},t;W^{1,p})} ),
\end{align*}
where $C$ is independent of $\wt{T}$ and $t$,
due to the global-in-time Strichartz estimates.
In view of \eqref{esti-wtvf.1}, \eqref{esti-wtvf.2} and \eqref{X'-X} below,
the first term on the right-hand side above is bounded by
\begin{align*}
    C \ve(\wt{T}) (\|\p_j \wt{u}\|_{LS{(\wt{T},t)}} + \| \wt{u}\|_{LS{(\wt{T},t)}} +  \| \wt{v}\|_{LS{(\wt{T},t)}}),
\end{align*}
which along with \eqref{h-0-Sigma}
and  the uniform bounds \eqref{glbdd-wtu-LpWq-Weight} and \eqref{glbdd-wtv-Sigma}
implies that $\bbp$-a.s.
\begin{align}  \label{ho-wtv**}
    \|\partial_j \wt{v}\|_{L^\g(\wt{T},t;L^\rho)} + \| \partial_j \wt{v}\|_{C([\wt{T},t];L^2)}
    \leq  C \ve(\wt{T}).
\end{align}

Thus, it follows from \eqref{ho-wtv*} and \eqref{ho-wtv**} that
\begin{align} \label{etsi-ve-wtv}
   \|\wt{v}\|_{L^\g(\wt{T},t;W^{1,\rho})} + \|\wt{v}\|_{C([\wt{T},t];H^1)}
    + \|\wt{v}\|_{LS{(\wt{T},t)}}
   \leq C \ve(\wt{T})
\end{align}
where $C$ is independent of $\wt{T}$ and $t$.
Letting $t\to 1^-$ we get
\begin{align} \label{asym-wtv-LpWq-Sigma}
    \|\wt{v}\|_{L^\g(\wt{T},1;W^{1,\rho})} + \|\wt{v}\|_{C([\wt{T},1);H^1)} + \|\wt{v}\|_{LS{(\wt{T},1)}}
    \leq  C \ve(\wt{T}) \to 0,\ \ as\ \wt{T}\to 1^{-},\ a.s..
\end{align}

Next, regarding the estimate of $\||\cdot| \wt{v}\|_{L^\g(0,1;L^\rho)}$,
similar to \eqref{equa-xwtz},
\begin{align*}
   \p_t (x_j \wt{v}) = \wt{A}_*(t) (x_j  \wt{v}) + H(\wt{v}, \p_j \wt{v}, \wt{u}),
\end{align*}
where $x_j\wt{v}(\wt{T}) =0$ and
\begin{align*}
   H(\wt{v}, \p_j \wt{v}, \wt{u})
   =& i (2 \p_j \wt{v} + \wt{b}_{*,j} \wt{v}) - i x_j (\wt{b}_* \cdot \na + \wt{c}_*) \wt{u} \\
    & -\lbb i h(t) (|e^{\wt{\vf}_*}(\wt{v}+\wt{u})|^{\a-1} x_j(\wt{v}+\wt{u}) - |\wt{u}|^{\a-1} (x_j \wt{u})).
\end{align*}
Then, similar  to \eqref{lpq-xwtz},
\begin{align} \label{esti-xwtv}
       & \|x_j \wt{v}\|_{L^\g(\wt{T},t ; L^\rho)}
         + \|x_j \wt{v}\|_{C([\wt{T},t];L^2)}  \nonumber \\
  \leq& C \bigg[  \|2 \p_j \wt{v} + \wt{b}_{*,j} \wt{v}\|_{L^1(\wt{T}, t; L^2)}
       + \|x_j (\wt{b}_* \cdot \na + \wt{c}_*) \wt{u}\|_{L^1(\wt{T}, t; L^2)} \nonumber  \\
      &\quad + \wt{\ve}_1(\wt{T}) \|\wt{v}+\wt{u}\|^{\a-1}_{L^q(\wt{T},t ; W^{1,p})}
          (\|x_j \wt{u}\|_{L^q(\wt{T},t ; L^p)} + \|x_j \wt{v}\|_{L^q(\wt{T},t ; L^p)} ) \nonumber  \\
      &\quad + \wt{\ve}_1(\wt{T}) \| \wt{u}\|^{\a-1}_{L^q(\wt{T},t ; W^{1,p})}
           \|x_j \wt{u}\|_{L^q(\wt{T},t ; L^p)}  \bigg], \ \ \bbp-a.s.,
\end{align}
where $C$ is independent of $\wt{T}$ and $t$.

Thus, by virtue of the uniform bounds \eqref{glbdd-wtu-LpWq-Weight}
and \eqref{glbdd-wtv-Sigma} we obtain
\begin{align} \label{asym-wtv-weight-Sigma}
    \|x_j \wt{v}\|_{L^\g(\wt{T},t ; L^\rho)}
         + \|x_j \wt{v}\|_{C([\wt{T},t];L^2)}
    \leq& C(1-\wt{T}) + C \wt{\ve}_1(\wt{T}) \to 0,\ \ as\ \wt{T}\to 1^{-},\ a.s..
\end{align}

Therefore,
combining \eqref{asym-wtv-LpWq-Sigma}
and \eqref{asym-wtv-weight-Sigma}
we obtain \eqref{asym-wtv-Sigma}.
The proof is finished.
\hfill $\square$

{\bf Proof of Theorem \ref{Thm-sca-X} (i).}
For any $t_1, t_2 $ close to $1$, we have
\begin{align*}
  |\wt{z}_*(t_1) - \wt{z}_*(t_2)|_{\Sigma}
  \leq |\wt{z}_*(t_1) - \wt{u}(t_1)|_{\Sigma}
       + |\wt{z}_*(t_2) - \wt{u}(t_2)|_{\Sigma}
       + |\wt{u}(t_1) - \wt{u}(t_2)|_{\Sigma},
\end{align*}
where $\wt{u}$ is the solution to \eqref{equa-wtu}
with $\wt{u}(\wt{T}) = \wt{z}_*(\wt{T})$,
$\wt{T}\in (0,1)$.

Moreover,
since by Proposition \ref{Prop-wtu},
$\wt{u}$ exists on $[0,1]$
and is continuous at time $1$,
we have $\bbp$-a.s.,
\begin{align*}
   \lim\limits_{t_1,t_2 \to 1^{-}}|\wt{u}(t_1) - \wt{u}(t_2)|_{\Sigma}
   =0.
\end{align*}

This implies that
\begin{align*}
  \limsup\limits_{t_1,t_2\to\9} |\wt{z}_*(t_1) - \wt{z}_*(t_2)|_{\Sigma}
  \leq& \limsup\limits_{t_1\to\9}   |\wt{z}_*(t_1) - \wt{u}(t_1)|_{\Sigma}
       + \limsup\limits_{t_2\to\9} |\wt{z}_*(t_2) - \wt{u}(t_2)|_{\Sigma}  \\
  \leq& 2  \|\wt{z}_* - \wt{u}\|_{C([\wt{T}, 1); \Sigma)}.
\end{align*}

Then, in view of \eqref{asym-wtv-Sigma},
letting $\wt{T} \to 1^-$ we obtain
\begin{align*}
   \lim\limits _{t_1,t_2\to\9} |\wt{z}_*(t_1) - \wt{z}_*(t_2)|_{\Sigma} =0.
\end{align*}
This yields that $\wt{z}_*$ has the limit at time $1$.
Thus, using the equivalence of asymptotics of $\wt{z}_*$  at time $1$
and $z_*$ at infinity
we conclude that $z_*$ scatters at infinity in the speudo-conformal space,
i.e., \eqref{sca-X-delta} holds.

Therefore,
the assertion of Theorem \ref{Thm-sca-X} $(i)$ is proved. \hfill $\square$ \\

Before proving Theorem \ref{Thm-sca-X} $(ii)$,
we obtain global estimates for the original solution $z_*$ to \eqref{equa-z*}
from those of its pseudo-conformal transformation $\wt{z}_*$.
\begin{corollary} \label{Cor-z*}
Let $\a\in (1+\a(d), 1+\frac{4}{d-2}]$, $d\geq 3$,
and $z_*$ be the solution to \eqref{equa-z*}.
Then, for any Strichartz pair $(\rho,\g)$,
\begin{align} \label{bdd-z-wpq}
     \|z_*\|_{L^\g(0,\9; W^{1,\rho})} <\9,\ \ a.s..
\end{align}
Moreover, if $\wt{q}:= \frac{2(\a^2-1)}{4-(\a-1)(d-2)}$, we have
\begin{align} \label{bdd-z-lap}
      \|z_*\|_{L^{\wt{q}}(0,\9;L^{\a+1})} \leq C < \9,\ \ a.s.,
\end{align}
where $C$ is independent of $\wt{T}$.
\end{corollary}

{\bf Proof.}
We shall prove the estimates of $z_*$ in \eqref{bdd-z-wpq} and \eqref{bdd-z-lap}
from those of the  pseudo-conformal transformation $\wt{z}_*$ in \eqref{glbdd-wtz-LpWq-Weight}.

For this purpose, since $(d/2-d/\rho)\g-2=0$,
straightforward computations show that
\begin{align} \label{bdd-z-wpq.1}
   \|z_*\|^\g_{L^\g(0,\9; L^{\rho})}
   = \int_0^1 (1-t)^{\frac{d(\rho-2)}{2\rho} \g -2} |\wt{z}_*(t)|^\g_{L^\rho} dt
   = \|\wt{z}_*\|^\g_{L^\g(0,1; L^{\rho})} <\9,\ \ a.s.,
\end{align}
where the last step is due to Corollary \ref{Cor-wtz-globH1}.
Moreover,  we compute that, if $s:= t/(1-t)$,
\begin{align*}
   |\na z_*(s)|_{L^\rho}
   = (1-t)^{\frac d2 +1 - \frac d \rho} |\na \wt{z}_*(t)- \frac i 2 \frac{x}{1-t} \wt{z}_*(t)|_{L^\rho}.
\end{align*}
Since
$(d/2 +1)\g - d\g /\rho -2 = \g$
and $d \g/2 - d \g /\rho -2 = 0$,
by Lemma \ref{Lem-h1-wtv},
this implies that
\begin{align} \label{bdd-z-wpq.2}
  \|\na z_*\|_{L^\g(0,\9;L^\rho)}
  \leq& C_\g\bigg( \int_0^1 (1-t)^{(\frac d2 +1)\g - \frac d \rho \g -2} |\na \wt{z}_*(t)|_{L^\rho}^\g dt \nonumber \\
      & \qquad +\int_0^1 (1-t)^{\frac d2 \g - \frac d \rho \g -2} ||\cdot|  \wt{z}_*(t)|^\g_{L^\rho} dt \bigg)^{\frac{1}{\g}} \nonumber \\
  \leq& C_\g(\|\na \wt{z}_*\|_{L^\g(0,1;L^\rho)} +  \||\cdot| \wt{z}_*\|_{L^\g(0,1;L^\rho)} )
   <\9,\ \  a.s..
\end{align}
Thus, taking together \eqref{bdd-z-wpq.1} and \eqref{bdd-z-wpq.2} yields \eqref{bdd-z-wpq}.

Regarding  \eqref{bdd-z-lap},
the argument is similar to that of \cite[Proposition $3.15 (iv)$]{CW92}.
Using Lemma \ref{Lem-h1-wtv} and the Sobolev imbedding theorem we have
\begin{align*}
     \sup \limits_{t\in[0,1]} |\wt{z}_*(t)|_{\a+1}
     \leq C \|\wt{z}_*\|_{C([0,1]; H^1)}
      < \9, \ \ a.s..
\end{align*}
Then, direct calculations imply that $\bbp$-a.s.
$$\|z_*\|^{\wt{q}}_{L^{\wt{q}}(0,\9;L^{\a+1})}
  = \int_0^\9 (1+s)^{-\frac{d(\a-1)}{2(\a+1)}\wt{q}} |\wt{z}_*(\frac{s}{1+s})|_{L^{\a+1}}^{\wt{q}} ds
  \leq C \int_0^\9 (1+s)^{- \frac{d(\a-1)}{2(\a+1)} \wt{q}} ds   <\9, $$
where the last step is due to the fact that
$-\frac{d(\a-1)}{2(\a+1)}\wt{q} = - \frac{d(\a-1)^2}{4-(\a-1)(d-2)}<-1$
for $\a\in (1+\a(d), 1+ 4/(d-2)]$.
Thus, we obtain \eqref{bdd-z-lap} and complete the proof.
\hfill $\square$ \\

{\bf Proof of Theorem \ref{Thm-sca-damped-X} (ii) (continued).}
Recall that $V(t,s)$, $s,t\geq 0$,
are evolution operators corresponding to \eqref{equa-z*} in the
homogeneous case $F\equiv 0$,
i.e., for any $v_s\in H^1$, $v(t):=V(t,s)v_s$ satisfies
\begin{align} \label{equa-homo-z*}
\p_t v(t) = A_*(t) v(t),\ \ v(s)= v_s.
\end{align}
We reformulate equation \eqref{equa-z*} in the mild form
\begin{align} \label{mild-wtz}
  z_*(t) = V(t,0) z_*(0) -\lbb i \int_0^t V(t,s) e^{-\vf_*(s)}F(e^{\vf_*(s)} z_*(s)) ds.
\end{align}
Then, for any $0<t_1<t_2<\9$,
if $w(t):= -\lbb i \int_{t_1}^{t} V(t,s) e^{- \vf_*(s)}F(e^{\vf_*(s)} z_*(s)) ds$,
we have
\begin{align} \label{diff-V}
  V(0,t_2) z_*(t_2) - V(0,t_1)  z_*(t_1)
  =&  V(0,t_2) w(t_2).
\end{align}

Moreover,
under the condition \eqref{asymflat}
and that $g_k\in L^2(\bbr^+)$, a.s.,
global-in-time Strichartz estimates hold for the operator $A_*$ in \eqref{equa-homo-z*}.
See Section \ref{Sec-Stri} below.

Below we consider two cases $\a\in (1+\a(d), 1+4/(d-2))$
and $\a= 1+4/(d-2)$ separately.

In the case where $\a\in (1+\a(d), 1+4/(d-2))$,
since
$V(\cdot,t_2) w(t_2)$ satisfies equation
\eqref{equa-homo-z*} with the final datum $w(t_2)$ at time $t_2$,
using \eqref{diff-V} and
applying  Corollary \ref{Cor-Stri} $(iii)$
with  $(p,q)= (\a+1, \frac{4(\a+1)}{d(\a-1)})$ we  obtain
\begin{align*}
  | V(0,t_2) z_*(t_2) - V(0,t_1)  z_*(t_1)  |_{H^1}
  \leq& \|V(\cdot,t_2) w(t_2) \|_{C([0,t_2];H^1)}
  \leq C | w(t_2) |_{H^1}, \ \ a.s.,
\end{align*}
where $C$ is independent of $t_1,t_2$.

Moreover, since $w(\cdot)$ satisfies \eqref{equa-z*} with
$w(t_1)=0$, and
$\vf_* \in L^\9(0,\9; W^{1,\9})$, a.s.,
applying Corollary \ref{Cor-Stri} $(ii)$  and H\"older's inequality we get
that the right-hand side above is bounded by
\begin{align} \label{V-t12}
    C \| w \|_{C([t_1,t_2];H^1)}
   \leq&   C \| e^{- \vf_*(s)}F(e^{\vf_*(s)} z_*) \|_{L^{q'}(t_1, t_2; W^{1,p'})} \nonumber \\
   \leq&   C\|z_*\|^{\a-1}_{L^{\wt{q}}(t_1,t_2;L^{p})}  \|z_*\|_{L^{q}(t_1,t_2; W^{1,p})},
\end{align}
where $\wt{q}= \frac{2(\a^2-1)}{4-(\a-1)(d-2)}$ and  $C$  is independent of $t_1, t_2$.

Thus, in view of \eqref{bdd-z-wpq} and \eqref{bdd-z-lap}, we obtain that
as $t_1, t_2 \to \9$,
\begin{align} \label{Cau-X}
     |V( 0,t_2)z_*( t_2) - V( 0,t_1)z_*(  t_1)|_{H^1} \to 0,\ \ a.s.,
\end{align}
which implies that $\bbp$-a.s. there exists $X_+ \in H^1$  such that
\begin{align}
     |V(  0,t)z_* (t) - X_+ |_{H^1} \to 0,\ \ as\ t\to \9,
\end{align}
thereby proving \eqref{sca-X-V} for the case where $\a\in (1+\a(d), 1+4/(d-2))$.

The energy-critical case where $\a=1+4/(d-2)$
can be proved similarly.
Choose the Strichartz pairs
$(p_1,p_1)= (2+4/d, 2+4/d)$
and $(p_2,q_2) = (\frac{2d(d+2)}{d^2+4}, \frac{2(d+2)}{d-2})$.
We only need to replace \eqref{V-t12} by the estimate below
\begin{align}
    | V(0,t_2) z_*(t_2) - V(0,t_1)  z_*(t_1)  |_{H^1}
   \leq&   C \|e^{- \vf_*(s)}F(e^{\vf_*(s)} z_*) \|_{L^{p_1'}(t_1, t_2; W^{1,p_1'})} \nonumber \\
   \leq&   C\|z_*\|^{\frac{4}{d-2}}_{L^{\frac{2(d+2)}{d-2}}((t_1,t_2)\times\bbr^d)}  \|z_*\|_{L^{p_1}(t_1,t_2; W^{1,p_1})} \nonumber \\
   \leq&   C\|z_*\|^{\frac{4}{d-2}}_{L^{q_2}(t_1,t_2; W^{1,p_2})}  \|z_*\|_{L^{p_1}(t_1,t_2; W^{1,p_1})},
\end{align}
where the last step is due to the Sobolev imbedding
$W^{1,\frac{2d(d+2)}{d^2+4}}  \hookrightarrow L^{\frac{2(d+2)}{d-2}}(\bbr^d)$.
Thus,
by virtue of \eqref{bdd-z-wpq}
we obtain \eqref{Cau-X},
thereby proving \eqref{sca-X-V} in the energy-critical case
where $\a=1+4/(d-2)$,
and the proof is complete. \hfill $\square$

\section{Scattering in the energy space} \label{Sec-H1}

In this section, we use the Strichartz pairs
$(p_1,p_1):= (2+\frac 4d,2+\frac 4d)$,
$(p_2,q_2):= ( \frac{2d(d+2)}{d^2+4},\frac{2(d+2)}{d-2} )$
and
define the spaces
$$Y^0(I):= L^{p_1} (I\times \bbr^d) \cap L^{q_2}(I; L^{p_2}),$$
and
$$Y^1(I) := \{f\in Y^0(I), \na f \in Y^0(I)\}.$$
Let  $p'_1 = \frac{2(d+2)}{d+4}$.
Then, $1/p_1' + 1/p_1 =1$.

\begin{lemma} \label{Lem-Inter-H1}
Let $\a \in [1+4/d, 1+4/(d-2)]$, $d\geq 3$,
and $I \times \bbr^d$ be an arbitrary spacetime slab.
Then, for any $f,g\in Y^1(I)$,
\begin{align}  \label{inter-H1}
        \||f|^{\a-1} g\|_{L^{p'_1}(I\times \bbr^d)}
 \leq  C \|f\|^{2- \frac{(\a-1)(d-2)}{2}}_{L^{p_1}(I\times \bbr^d)}
        \|f\|^{\frac{d(\a-1)}{2}-2}_{L^{q_2}(I; W^{1, p_2})}
        \| g\|_{L^{p_1}(I\times \bbr^d)}.
\end{align}
Moreover, if in addition $\a \geq 2$, we have
\begin{align}  \label{inter-H1*}
        \||f|^{\a-2} \na f g\|_{L^{p'_1}(I\times \bbr^d)}
 \leq  C \|f\|^{\a-1}_{Y^1(I)}
         \| g\|_{Y^1(I)}.
\end{align}
\end{lemma}

{\bf Proof.}
By Lemma 2.6 in \cite{TVZ07} we know that,
if $\theta_1:= \frac{2}{\a-1}- \frac{d-2}{2},
\theta_2:=\frac{d}{2}-\frac{2}{\a-1} \in [0,1]$,
\begin{align*}
    \||f|^{\a-1}  g\|_{L^{\frac{2(d+2)}{d+4}}(I\times \bbr^d)}
   \leq C \|f\|^{(\a-1)\theta_1}_{L^{2+\frac 4d}(I\times \bbr^d)}
        \|f\|^{(\a-1)\theta_2}_{L^{\frac{2(d+2)}{d-2}}(I\times \bbr^d)}
        \| g\|_{L^{2+\frac 4d}(I\times \bbr^d)},
\end{align*}
which along with the Sobolev imbedding
$W^{1, \frac{2d(d+2)}{d^2+4}} \hookrightarrow  L^{\frac{2(d+2)}{d-2}}$ implies \eqref{inter-H1}.

Regarding \eqref{inter-H1*},
using H\"older's inequality,
we have
\begin{align*}
        \||f|^{\a-2} \na f g\|_{L^{p'_1}(I\times \bbr^d)}
        \leq&   \|f\|^{\a-2}_{L^\frac{(d+2)(\a-1)}{2}(I\times \bbr^d)}
              \| \na f\|_{L^{2+\frac4d}(I\times \bbr^d)}
              \|g\|_{L^\frac{(d+2)(\a-1)}{2}(I\times \bbr^d)}.
\end{align*}
Note that,
\begin{align*}
    \|f\|_{L^\frac{(d+2)(\a-1)}{2}(I\times \bbr^d)}
    \leq \|f\|^{\theta_1}_{L^{2+\frac 4d}(I\times \bbr^d)}
        \|f\|^{\theta_2}_{L^{\frac{2(d+2)}{d-2}}(I\times \bbr^d)},
\end{align*}
and similar estimate holds also for $g$.
We obtain
\begin{align*}
    \||f|^{\a-2} \na f g\|_{L^{p'_1}(I\times \bbr^d)}
       \leq&  \|f\|^{(\a-2)\theta_1}_{L^{2+\frac 4d}(I\times \bbr^d)}
              \|f\|^{(\a-2)\theta_2}_{L^{\frac{2(d+2)}{d-2}}(I\times \bbr^d)}
              \| \na f\|_{L^{2+\frac4d}(I\times \bbr^d)}\\
       &    \cdot
              \|g\|^{\theta_1}_{L^{2+\frac 4d}(I\times \bbr^d)}
              \|g\|^{\theta_2}_{L^{\frac{2(d+2)}{d-2}}(I\times \bbr^d)},
\end{align*}
which implies \eqref{inter-H1*},
thereby finishing the proof.
\hfill $\square$ \\

\begin{lemma} \label{Lem-Ito-L2H1}
Consider the situation of Theorem \ref{Thm-Sca-H1}
and let $z_*$ be the unique global solution  to \eqref{equa-z*}
with $z_*(0)=X_0$.
Then, for each $X_0\in H^1$, we have
\begin{align}  \label{glo-z*-H1}
   \sup\limits_{0\leq t<\9} |z_*(t)|_{H^1} \leq C <\9,\ \ a.s..
\end{align}
\end{lemma}

{\bf Proof.} \eqref{glo-z*-H1} follows directly from
\eqref{glo-X-H1} and the fact that $\vf_*\in L^\9(\bbr^+;W^{1,\9})$, a.s..
\hfill $\square$

\begin{remark}
Unlike in the pseudo-conformal case,
the uniform bound \eqref{glo-z*-H1} of $z_*$
does not imply immediately the scattering behavior of $z_*$
in the energy space.
Actually, in the deterministic case,
more subtle interaction Morawetz estimates
are needed to obtain the scattering.
See, e.g., \cite{CKSTT08,RV07,V07}.
\end{remark}

Below, we prove the  scattering
by using the idea of comparison as in the proof of Lemma \ref{Lem-h1-wtv}.
For this purpose, we first prove the uniform global estimates
(independent of $T$) for the solution $u$ to \eqref{equa-u}.
\begin{proposition} \label{Pro-Sca-Det}
Consider $\lbb =-1$,  $\a\in [1+4/d,1+4/(d-2)]$, $d\geq 3$.
Then, if $u(T)= z_*(T) \in H^1$,
there exists a unique global $H^1$-solution $u$ (depending on $T$) to \eqref{equa-u} such that
$u$ scatters at infinity in $H^1$ and
for any $1\leq j\leq d$ and any Strichartz pair $(\rho, \g)$,
\begin{align} \label{gloStri-u-H1}
    \|u\|_{L^\g(\bbr^+;W^{1,\rho})} + \|u\|_{LS(\bbr^+)} + \|\partial_j u\|_{LS(\bbr^+)} \leq C<\9,
\end{align}
where $LS(\bbr^+)$ is the local smoothing space defined in Section \ref{Sec-Stri},
and $C$ is independent of $T$.
\end{proposition}

{\bf Proof. }
We mainly consider the global well-posedness, scattering and global bound of $ \|u\|_{L^\g(\bbr;W^{1,\rho})}$.
The global bound in the local smoothing space can be proved standardly by Strichartz estimates.

The crucial estimates used below are that,
for $\a\in (1+4/d, 1+4/(d-2)]$ and any Strichartz pair $(\rho, \g)$,
\begin{align} \label{u-LpWq-L2H-H1}
      \|u\|_{L^\g(\bbr^+; W^{1,\rho})}
    \leq C(|u(T)|^2_2, H(u(T))),
\end{align}
where $H(u)$ is the Hamiltonian of $u$,
i.e., $H(u) = \frac 12 |\na u|_2^2 + \frac {1}{\a+1} |u|_{L^{\a+1}}^{\a+1}$,
and for the mass-critical case $\a=1+ 4/d$,
if $p :=2+4/d$,
\begin{align} \label{u-Lp-L2-L2}
   \|u\|_{L^{p}(\bbr^+\times\bbr^d)} \leq C(|u(T)|^2_2).
\end{align}
See \cite[Subsection 5.3]{TVZ07} for the inter-critical case
$\a\in (1+4/d, 1+4/(d-2))$,
\cite{CKSTT08,RV07,V07} and \cite[Subsection 4.2]{TVZ07}
for the energy-critical case,
and \cite{D12} for the mass-critical case.

Now, in the case where $\a\in (1+4/d, 1+4/(d-2)]$,
if $u(T) = z_*(T) \in H^1$,
the global well-posedness and scattering of $u$
(with $T$ fixed)
follow from \cite{C03,CKSTT08,RV07,V07}.
Moreover, by virtue of \eqref{u-LpWq-L2H-H1} and \eqref{glo-z*-H1},
we have for any Strichartz pair $(\rho, \g)$, $\bbp$-a.s.,
\begin{align} \label{LqWp-L2H-H1}
    \|u\|_{L^\g(\bbr^+; W^{1,\rho})}
    \leq C(|u(T)|^2_2, H(u(T)))
    = C(|z_*(T)|^2_2, H(z_*(T)))
    \leq C(|z_*(T)|_{H^1} )
    \leq C
\end{align}
where $C$ is independent of $T$.

Regarding the mass-critical case where $\a =1+ 4/d$,
since $u(T) = z_*(T) \in H^1 \subset L^2$,
it follows from \cite{D12} that
there exists a unique global $L^2$-solution to \eqref{equa-u}.
Moreover, by \eqref{u-Lp-L2-L2} and \eqref{glo-z*-H1},
\begin{align} \label{u-LpLp}
   \|u\|_{L^{p}(\bbr^+\times\bbr^d)} \leq C(|u(T)|^2_2) \leq C(|z_*(T)|^2_2) \leq C <\9,\ \ a.s.,
\end{align}
where $C$ is independent of $T$.
Then,
we split $\bbr^+$ into $\bigcup_{j=0}^M(t_j,t_{j+1})$ with
$t_0 = 0$, $t_{M+1}=\9$,
$M (= M(T))<\9$,
such that $\|u\|_{L^{p}((t_j,t_{j+1})\times\bbr^d)}=\ve$,
$0\leq j\leq M-1$,
and
$\|u\|_{L^{p}((t_M,\9)\times\bbr^d)} \leq \ve$,
where $\ve$ is to be chosen later.
By \eqref{u-LpLp},
\begin{align} \label{glbdd-MT}
   \sup\limits_{0<T<\9} M(T) \leq (\frac C \ve)^{p}
\end{align}
with $C$ independent of $T$.
Note that,
the conservations of the Hamiltonian
$H(u(T)) = H(u(t))$ and
the mass $|u(T)|_2^2 = |u(t)|_2^2$, $t\in \bbr$,
together with \eqref{glo-z*-H1}  imply that,
for $0<C<\9$ independent of $T$,
\begin{align} \label{na-Ham}
|u(t)|_{H^1} \leq   |u(t)|_2^2 + 2H(u(t))
             = |z_*(T)|_2^2 + 2H(z_*(T))
              \leq C (1+\|z_*\|^{\a+1}_{C(\bbr^+; H^1)}) <\9,\ \ a.s..
\end{align}
Then,
applying Strichartz estimates to \eqref{equa-u} and using H\"older's inequality
we get that, if $p':=2(d+2)/(d+4)$, then
\begin{align} \label{gloStri-u-H1*}
   \|u\|_{L^{p}(t_j,t_{j+1};W^{1,p})}
   \leq& C |u(t_j)|_{H^1} + C \||u|^{\a-1}u\|_{L^{p'}(t_j,t_{j+1}; W^{1,p'})} \nonumber  \\
   \leq& C |u(t_j)|_{H^1}
         + C \|u\|^{\frac 4 d}_{L^{p}(t_j,t_{j+1}; L^{p})} \|u\|_{L^{p}(t_j,t_{j+1}; W^{1,p})} \nonumber \\
   \leq& C + C \ve^{\frac 4 d} \|u\|_{L^{p}(t_j,t_{j+1}; W^{1,p})},\ \ \forall\ 0\leq j\leq M.
\end{align}
Thus, taking $\ve$ small enough (independent of $T$), summing over $j$
and using the global bound \eqref{glbdd-MT}
we obtain that $ \|u\|_{L^{p}(\bbr^+;W^{1,p})} \leq C(M) \leq C<\9$,
where $C$ is independent of $T$.
Hence, $ \|u\|_{L^\g(\bbr^+;W^{1,\rho})} \leq C$
for any Strichartz pair $(\rho, \g)$ by Strichartz estimates.

In particular,
this yields the scattering of $u$ in $H^1$.
Actually, by \eqref{equa-u},
\begin{align*}
   e^{it_2 \Delta}u(t_2) - e^{it_1 \Delta}u(t_1)
   = -\lbb i \int_{t_1}^{t_2} e^{is\Delta} F(u(s))ds, \ \ \forall \ t_1,t_2\geq T>0.
\end{align*}
Similar  to \eqref{gloStri-u-H1*}, we have as $t_1,t_2 \to \9$,
\begin{align} \label{diff-e-u}
   |e^{it_2 \Delta}u(t_2) - e^{it_1 \Delta}u(t_1)|_{H^1}
   \leq C \|u\|^{\frac{4}{d}}_{L^{p}(t_1,t_2; L^{p})} \|u\|_{L^{p}(t_1,t_2; W^{1,p})} \to 0,
\end{align}
which implies the scattering in $H^1$ for the mass-critical exponent $\a=1+4/d$.

Therefore, the proof is complete.
\hfill $\square$

For the solution $v$ to \eqref{equa-v},
we have the crucial global estimates and asymptotics at infinity below.
\begin{lemma} \label{Bdd-LqW1p}
Assume the conditions of Theorem \ref{Thm-Sca-H1} to hold.
Let $v $ be the solution to \eqref{equa-v}
with $v(T) =0$.
Then, $\bbp$-a.s., for any Strichartz pair $(\rho,\g)$,
\begin{align} \label{glbdd-v-LqWp-H1}
     \|v\|_{L^\g(T,\9; W^{1,\rho})} \leq C <\9,\
\end{align}
where $C$ is uniformly bounded for $T$ large enough,
and  for any Strichartz pair $(\rho, \g)$,
as $T\to \9$,
\begin{align} \label{v-0-H1}
   \|v\|_{L^\g(T,\9; W^{1,\rho})} + \|v\|_{C([T,\9); H^1)} \to 0.
\end{align}
\end{lemma}

We remark that,
unlike in the pseudo-conformal case,
the global estimate \eqref{glo-z*-H1} dose not imply directly
the global bound of $z_*$ (and also $v$) in
the Strichartz space $L^\g(T,\9; W^{1,\rho})$.
Proceeding differently,
we will prove \eqref{glbdd-v-LqWp-H1}
by using the comparison arguments
which  simultaneously gives \eqref{v-0-H1} as well.

{\bf Proof.}
As mentioned in the proof of Theorem \ref{Thm-sca-X} $(ii)$,
global-in-time Strichartz estimates
hold for the operator $A_*$ in \eqref{equa-v},
due to the fact that  $g_k\in L^2(\bbr^+)$, $\bbp$-a.s., $1\leq k\leq N$ (see Section \ref{Sec-Stri} below).

Below,
let $(p_i,q_i)$, $i=1,2$,
and $Y^1(I)$, $I \subset \bbr^+$, be the Strichartz pairs and Strichartz space, respectively,
as in Lemma \ref{Lem-Inter-H1}.

Applying Corollary \ref{Cor-Stri} to \eqref{equa-v}
we have  for any  $0<T<t<\9$,
\begin{align} \label{esti-v-H1}
      & \|v\|_{C([T, t]; H^1)} + \|v\|_{Y^1(T,t)}  \nonumber \\
  \leq& C \| (b_*\cdot \na +c_*)u \|_{LS'{(T,t)}}
        + C \| \na ((b_*\cdot \na +c_*) u) \|_{LS'{(T,t)}} \nonumber \\
      & + C\| (e^{(\a-1)\Re \vf_*}-1) F(u+v) \|_{L^{p_1'}(T,t; W^{1,p_1'})} \nonumber \\
      & + C\| F(u+v) - F(u) \|_{L^{p_1'}(T,t; W^{1,p_1'})}
\end{align}
where $C$ is independent of $T$ and $t$,
due to the global-in-time Strichartz estimates.

Below,
we estimate the four terms on the right-hand side of \eqref{esti-v-H1} respectively.

First, since $g_k\in L^\9(\Omega; L^2(\bbr^+))$,
using \eqref{asymflat} we have (see Section \ref{Sec-Stri} below) that $\bbp$-a.s.
for $0\leq |\beta|\leq 2$,  $0\leq |\g|\leq 1$,
\begin{align}  \label{ve-H1}
\ve(T) := \sup_{t\in [T,\9),x\in\bbr^d} \<x\>^2(|\partial_x^\beta b_*(t,x)|+|\partial_x^\g c_*(t,x)|) \to 0, \ as\ T\to \9,
\end{align}
and
$$ \limsup_{|x|\to \9}\sup\limits_{0\leq t<\9} \<x\>^2(|\partial_x^\beta b_*(t,x)|+|\partial_x^\g c_*(t,x)|) = 0,  $$
which implies by \eqref{X'-X} below (see also Remark \ref{Rem-Stri})
and Proposition \ref{Pro-Sca-Det} that
\begin{align} \label{esti-v-H1.1}
   &\| (b_*\cdot \na +c_*)u \|_{LS'{(T,t)}}
        +  \| \na ((b_*\cdot \na +c_*)u) \|_{LS'{(T,t)}}  \nonumber \\
   \le& C \ve(T) ( \|u\|_{LS{(T,t)}} +  \|\na u\|_{LS{(T,t)}})
   \leq C \ve(T),
\end{align}
where $C$ is independent of $T$.

Moreover, using the inequality $|e^x-1|\leq e|x|$ for $|x|\leq 1$
and the fact that $\|\vf_*\|_{C([T,\9);W^{1,\9})} \leq C \ve(T)$
we get
\begin{align} \label{esti-v-H1.2}
    \| (e^{(\a-1) \Re \vf_*}-1) F(u+v) \|_{L^{p_1'}(T,t; W^{1,p_1'})}
    \leq& C \ve(T) (\|u\|^\a_{Y^1(T,t)} + \|v\|^\a_{Y^1(T,t)} )  \nonumber \\
    \leq& C \ve(T) (1+ \|v\|^\a_{Y^1(T,t)} ),
\end{align}
where $C$ is independent of $T$,
due to Proposition \ref{Pro-Sca-Det}.

Regarding the fourth term on the right-hand side of \eqref{esti-v-H1},
we first note that
\begin{align} \label{esti-F-1}
   |F(u+v)- F(u)| \leq C (|u|^{\a-1} + |v|^{\a-1}) |v|.
\end{align}
Moreover, since
$\na F(u) = F_z(u) \na u + F_{\ol{z}}(u) \na \ol{u}$,
where $F_z, F_{\ol{z}}$ are the usual complex derivatives,
and similar equality holds for $\na F(u+v)$,
we have
\begin{align*}
    \na( F(u+v) - F(u))
    =& (F_z(u+v) - F_z(u)) \na u
      + (F_{\ol{z}}(u+v) - F_{\ol{z}}(u)) \na \ol{u} \\
    &  + F_z(u+v) \na v + F_{\ol{z}}(u+v) \na \ol{v}.
\end{align*}
Note that
$|F_z(u+v)| + |F_{\ol{z}}(u+v)| \leq C |u+v|^{\a-1}$,
and since $\a\geq 2$,
$|F_z(u+v) - F_z(u)| + |F_{\ol{z}}(u+v) - F_{\ol{z}}(u)|
  \leq C(|u|^{\a-2} + |v|^{\a-2}) |v|$.
Hence,
\begin{align} \label{esti-F-2}
    |\na( F(u+v) - F(u))|
    \leq (|u|^{\a-2} + |v|^{\a-2}) |v| |\na u|
        + (|u|^{\a-1}+|v|^{\a-1}) |\na v|,
\end{align}
Then, using \eqref{esti-F-1}, \eqref{esti-F-2} and Lemma \ref{Lem-Inter-H1} we obtain
\begin{align} \label{esti-v-H1.3}
   \| F(u+v) - F(u) \|_{L^{p_1'}(T,t; W^{1,p_1'})}
   \leq& C(\|u\|_{Y^1(T,t)} \|v\|^{\a-1}_{Y^1(T,t)}
           +\|u\|^{\a-1}_{Y^1(T,t)}  \|v\|_{Y^1(T,t)} \nonumber \\
       &\ \     + \|v\|^\a_{Y^1(T,t)}),
\end{align}
where $C$ is independent of $T$ and $t$.

Thus, plugging \eqref{esti-v-H1.1}, \eqref{esti-v-H1.2} and \eqref{esti-v-H1.3} into \eqref{esti-v-H1} we obtain
\begin{align} \label{esti-v-H-0}
   \|v\|_{C([T, t]; H^1)} + \|v\|_{Y^1(T,t)}
   \leq& C \ve(T) + C \|v\|^\a_{Y^1(T,t)} \nonumber \\
       & +C(\|u\|_{Y^1(T,t)} \|v\|^{\a-1}_{Y^1(T,t)}
           +\|u\|^{\a-1}_{Y^1(T,t)}  \|v\|_{Y^1(T,t)}).
\end{align}

Below we shall prove the   estimate
\begin{align} \label{glo-v-H1}
    \|v\|_{C([T,\9);H^1)} + \|v\|_{Y^1(T,\9)} \leq C' \ve(T),
\end{align}
where $C'$ is independent of $T$.

For this purpose,
we choose $\eta, \ve(T)>0$ sufficiently small, such that
\begin{align} \label{cond-eta-ve}
    D_0(\eta + \eta^{\a-1}) \leq \frac 12, \ \
    4D_0 D^*(\eta) \ve(T) \leq (1-\frac 1\a) (2\a D_0)^{-\frac{1}{\a-1}},
\end{align}
where the constants $D_0, D^*(\eta)$
are as in \eqref{esti-v-H-1} and \eqref{D*-H1}
below respectively, independent of $T$.

For each $T$ fixed,
since $\|u\|_{Y^1(T,\9)} \leq C_* <\9$,
we can choose $\{t_{1,j}\}_{j=0}^{M_1+1}$,
$t_{1,0} = T$, $t_{1,M_1+1}=\9$,  $M_1 (=M_1(T))<\9$,
such that $\|u\|_{L^{p_1}(t_{1,j}, t_{1,j+1})} = \frac \eta 2$,
$0\leq j\leq M_1-1$,
and $\|u\|_{L^{p_1}(t_{1,M_1}, \9)} \leq \frac \eta 2$.
Then, similarly to \eqref{glbdd-MT},
$\sup_T M_1(T) \leq (\frac{2C_*}{\eta})^{p_1}$.
Similarly,
we choose another partition $\{t_{2,j}\}_{j=0}^{M_2}$,
such that $\|u\|_{L^{q_2}(t_{2,j}, t_{2,j+1}; W^{1,p_2})} = \frac \eta 2$,
$0\leq j\leq M_2-1$,
$\|u\|_{L^{q_2}(t_{2,M_2}, \9; W^{1,p_2})} \leq \frac \eta 2$,
and so
$\sup_T M_2(T) \leq (\frac{2C_*}{\eta})^{q_2}$.

Thus, we can divide $\bbr^+$ into
finitely many subintervals $\{[t_j, t_{j+1}]\}_{j=0}^M$,
satisfying that
$\{t_j\}_{j=0}^{M+1} = \{t_{1,j}\}_{j=0}^{M_1+1} \cup \{t_{2,j}\}_{j=0}^{M_2+1}$,
$\|u\|_{Y^1(t_j, t_{j+1})} \leq \eta$,
$0\leq j\leq M$,
and
\begin{align} \label{bdd-MT-H1}
     \sup\limits_{0<T<\9} M(T)
     \leq \sup\limits_{0<T<\9}(M_1(T) + M_2(T))
     \leq (\frac{2C_*}{\eta})^{p_1} + (\frac{2C_*}{\eta})^{q_2}
     := C^*(\eta)<\9.
\end{align}

Now, let us  consider the estimate on the time interval $(T,t_1)$.
By \eqref{esti-v-H-0}
and the fact that $\a \geq 2$,
\begin{align} \label{esti-v-H-1}
   &\|v\|_{C([T, t_1]; H^1)} + \|v\|_{Y^1(T,t_1)} \nonumber \\
   \leq& C \ve(T) + C \|v\|^\a_{Y^1(T,t_1)}
        +C \eta \|v\|^{\a-1}_{Y^1(T,t_1)}
           + C \eta^{\a-1} \|v\|_{Y^1(T,t_1)} \nonumber \\
   \leq& D_0 \ve(T) + D_0  (\eta + \eta^{\a-1}) \|v\|_{Y^1(T,t_1)}
          + D_0 \|v\|^\a_{Y^1(T,t_1)},
\end{align}
where $D_0 = 2 C (\geq 1) $ is independent of $T$.
Then, by the choice of $\eta$ in \eqref{cond-eta-ve},
\begin{align*}
   \|v\|_{C([T, t_1]; H^1)} + \|v\|_{Y^1(T,t_1)}
   \leq  2D_0 \ve(T)  +  2D_0  \|v\|^\a_{Y^1(T,t_1)}.
\end{align*}
This, via \cite[Lemma A.1]{BRZ16.2} and
implies that
for $\ve(T)$ small enough such that
$2D_0\ve(T) < (1-\frac 1\a) (2\a D_0)^{-\frac{1}{\a-1}}$,
\begin{align} \label{esti-v-Tt1}
   \|v\|_{C([T, t_1]; H^1)} + \|v\|_{Y^1(T,t_1)}
   \leq  D_1 \ve(T)
\end{align}
with $D_1:= 2(4D_0)^{1+\a} (\frac{\a}{\a-1})^\a$, independent of $T$.
In particular,
$|v(t_1)|_{H^1} \leq D_1 \ve(T)$.

Next we use the inductive arguments.
Suppose that at the $j$-th step
$|v(t_j)|_{H^1} \leq D_j \ve(T)$, $0<j\leq M$,
where $D_j (\geq 1)$ is increasing with $j$
and is independent of $T$.

Then, similarly to \eqref{esti-v-H-1}, we have
\begin{align*}
   &\|v\|_{C([t_j, t_{j+1}]; H^1)} + \|v\|_{Y^1(t_j, t_{j+1})} \nonumber \\
   \leq& C|v(t_j)|_{H^1}
        + D_0 \ve(T) + D_0 (\eta + \eta^{\a-1}) \|v\|_{Y^1(t_j, t_{j+1})}
          + D_0 \|v\|^\a_{Y^1(t_j, t_{j+1})}.
\end{align*}
By the choice of $\eta$ in \eqref{cond-eta-ve} and the inductive assumption, we get
\begin{align*}
   \|v\|_{C([t_j, t_{j+1}]; H^1)} + \|v\|_{Y^1(t_j, t_{j+1})}
   \leq& 2 C D_j \ve(T)
        + 2D_0 \ve(T) +  2D_0 \|v\|^\a_{Y^1(t_j, t_{j+1})} \nonumber \\
   \leq& 4 D_0D_j \ve(T) + 2D_0 \|v\|^\a_{Y^1(t_j, t_{j+1})}.
\end{align*}
Using \eqref{cond-eta-ve} and apply \cite[Lemma A.1]{BRZ16.2} again,
we have that
for $\ve(T)$ even smaller such that
$4D_0 D_j \ve(T) < (1-\frac1\a)(2\a D_0)^{-\frac{1}{\a-1}}$,
\begin{align} \label{esti-v-H-j1}
   \|v\|_{C([t_j, t_{j+1}]; H^1)} + \|v\|_{Y^1(t_j, t_{j+1})}
   \leq  D_{j+1} \ve(T)
\end{align}
with $D_{j+1} = 2(4D_0)^{1+\a}(\frac{\a}{\a-1} D_j)^\a = D_1 D_j^\a$ and
is  independent of $T$.
In particular,
$|v(t_{j+1})|_{H^1} \leq D_{j+1} \ve(T)$.

Thus, letting
\begin{align} \label{D*-H1}
    D^*(\eta)= D_1^{\frac{\a^{C^*(\eta)+1}-1}{\a-1}} < \9
\end{align}
with $C^*(\eta)$ as in \eqref{bdd-MT-H1} above,
we deduce from the condition \eqref{cond-eta-ve}
and inductive arguments that,
for each $0\leq j\leq M \leq C^*(\eta)$,
the estimate \eqref{esti-v-H-j1} is valid
and
\begin{align*}
    D_{j+1} = D_1^{1+\a+\cdots +\a^{j}}
            = D_1^{\frac{\a^{j+1}-1}{\a-1}}
            \leq D^*(\eta) <\9,
\end{align*}
provided $\ve(T)$ satisfies the smallness condition in \eqref{cond-eta-ve} above.

Therefore, taking the sum or maximum over $0\leq j\leq M$
and taking into account
the uniform bound \eqref{bdd-MT-H1}
we obtain  \eqref{glo-v-H1}, as claimed.

Finally, \eqref{glbdd-v-LqWp-H1} and \eqref{v-0-H1} follow from \eqref{ve-H1}, \eqref{glo-v-H1}
and Strichartz estimates. The proof is complete.
 \hfill $\square$ \\

{\bf Proof of Theorem \ref{Thm-Sca-H1}.}
$(i)$.
Note that,  for any $t_1,t_2\geq T$
\begin{align*}
     | e^{it_1\Delta} z_*(t_1) - e^{it_2 \Delta} z_*(t_2) |_{H^1}
    \leq& | e^{it_1\Delta}( z_*(t_1) - u(t_1) ) |_{H^1}
        + | e^{it_2\Delta} (u(t_2) - z_*(t_2)) |_{H^1}   \\
        & + | e^{it_1\Delta} u(t_1) - e^{it_2 \Delta} u(t_2) |_{H^1} \\
    \leq& 2 \|z_* - u\|_{C([T,\9);H^1)}
           + | e^{it_1\Delta} u(t_1) - e^{it_2 \Delta} u(t_2) |_{H^1},
\end{align*}
where $u$ is the solution to \eqref{equa-u} with
$u(T) = z_*(T)$.

Since by Proposition \ref{Pro-Sca-Det},
for each $T$ fixed,
$u$ scatters
at infinity,
we have
$$ \lim\limits_{t_1,t_2\to \9} | e^{it_1\Delta} u(t_1) - e^{it_2 \Delta} u(t_2) |_{H^1}  = 0. $$

Then,  we get
\begin{align*}
   \limsup\limits_{t_1,t_2\to \9} | e^{it_1\Delta} z_*(t_1) - e^{it_2 \Delta} z_*(t_2) |_{H^1}
   \leq 2 \|z_* - u\|_{C([T,\9);H^1)}.
\end{align*}

Thus,
by virtue of \eqref{v-0-H1}, we taking $T$ to infinity to obtain that
$\{e^{it\Delta} z_*(t)\}$ is a Cauchy sequence in the space $H^1$,
which implies the scattering of $z_*$ at infinity specified in \eqref{sca-H1}.

$(ii)$. As in the proof of Theorem \ref{Thm-sca-X} $(ii)$,
Strichartz estimates and Lemma \ref{Lem-Inter-H1} imply  that $\bbp$-a.s.,
\begin{align*}
      &| V(0,t_2) z_*(t_2) - V(0,t_1)  z_*(t_1)  |_{H^1} \nonumber \\
  \leq&  C   \|e^{- \vf_*(s)}F(e^{\vf_*(s)} z_*)   \|_{L^{\frac{2(d+2)}{d+4}}(t_1, t_2; W^{1,\frac{2(d+2)}{d+4}})}
   \leq C\|z_*\|^{\a}_{Y^1(t_1,t_2)} ,
\end{align*}
where $Y^1(t_1,t_2)$ is the space defined in the previous proof of Lemma \ref{Bdd-LqW1p}
and $C$ is independent of $t_1,t_2$,
due to the global-in-time Strichartz estimates for $A_*$.

Note that,
$z_*= e^{-\vf_*}X= v+u$.
We have
$\|z_*\|_{Y^1(T,\9)} <\9$
for $T$ large enough,
due to the global bounds \eqref{gloStri-u-H1} and \eqref{glbdd-v-LqWp-H1}.

Thus, it follows that
$| V(0,t_2) z_*(t_2) - V(0,t_1)  z_*(t_1)  |_{H^1} \to 0$,
as $t_1,t_2\to \9$, $\bbp$-a.s.,
thereby yielding \eqref{sca-X-VH1}.
Therefore, the proof is complete. \hfill $\square$

\section{Proof of Theorem \ref{Thm-sca-damped-X}.} \label{Sec-Sca-Damp}

We first consider the scattering in the pseudo-conformal space.
As in Section \ref{Sec-Finite-Varia},
let $\wt{z}$ be the pseudo-conformal transformation of the
solution $z$ to \eqref{equa-z}, i.e.,
\begin{align}
   \wt{z}(t,x) := (1-t)^{-\frac d2}\ z\(\frac{t}{1-t}, \frac{x}{1-t}\)\ e^{i\frac{|x|^2}{4(1-t)}},
\end{align}
where $t\in [0,1)$, $x\in \bbr^d$.
We have
\begin{align} \label{equa-wtz-damp}
   &\p_t \wt{z} = - i \Delta \wt{z} - \lbb i h(t) e^{-\wt{\vf}} F( e^{\wt{\vf}} \wt{z} ),  \\
   & \wt{z}(0) = X_0 e^{i\frac{|x|^2}{4}}\in \Sigma, \nonumber
\end{align}
where $h(t)$ is as in Section \ref{Sec-Finite-Varia},
i.e., $h(t)= (1-t)^{\frac{d(\a-1)-4}{2}}$,
and
$$\wt{\vf}(t) = \vf(\frac{t}{1-t})
  =\sum\limits_{k=1}^N \(\int_0^{\frac{t}{1-t}} v_k g_k(s) d\beta_k(s) -  \int_0^{\frac{t}{1-t}} (\Re v_k) v_k g^2_k(s) ds\).$$

We  show that, for $\Re v_1$ large enough,
$\wt{z} (= \wt{z}_{v_1}) $ exists on $[0,1]$ with high probability.

For this purpose,
set $\wt{\calx}_\sigma^M := \{w \in L^\9(0,\sigma; L^2) \cap L^q(0,\sigma; L^p): \|w\|_{L^\9(0,\sigma; H^1)} + \|w\|_{L^q(0,\sigma; W^{1,p})} \leq M \}$
with $(p,q)=(\frac{d(\alpha+1)}{d+\alpha-1},\frac{4(\alpha+1)}{(d-2)(\alpha-1)})$,
and define $\wt{\Phi}$ on $\wt{\calx}_\sigma^M$  by
\begin{align*}
    \wt{\Phi}(w)(t) = e^{-it \Delta} (X_0 e^{i\frac{|x|^2}{4}}) - \lbb i \int_0^t e^{-i(t-s) \Delta} h(s) e^{-\wt{\vf}(s)} F(e^{\wt{\vf}(s)} w(s)) ds,
     \ \ w\in \wt{\calx}_\sigma^M.
\end{align*}
Then, similarly to \eqref{z*-lqp},  for any $w_j \in \wt{\calx}_\sigma^M$, $j=1,2$,
\begin{align} \label{Phi-v.1}
      \|\wt{\Phi}(w_j)\|_{L^\9(0,\sigma; H^1)} + \|\wt{\Phi}(w_j)\|_{L^q(0,\sigma; W^{1,p})}
      \leq   C|X_0|_{\Sigma} +  C \wt{\ve}^{\frac 1 \theta}_\sigma(v_1) M^\a,
\end{align}
and
\begin{align}\label{Phi-v.2}
      &\|\wt{\Phi}(w_1) - \wt{\Phi}(w_2)\|_{L^\9(0,\sigma; L^2)}
      +  \|\wt{\Phi}(w_1) - \wt{\Phi}(w_2)\|_{L^q(0,\sigma; L^p)} \nonumber \\
      \leq& C \wt{\ve}_\sigma^{\frac 1 \theta}(v_1)   M^{\a-1} \|w_1 -w_2\|_{L^q(0,\sigma; L^p)}.
\end{align}
where $\wt{\ve}_\sigma(v_1):= |he^{(\a-1) \Re  \wt{\vf}} |^\theta_{L^\theta(0,\sigma)}$,
$1<\theta<\9$ is as in \eqref{z*-lqp},
and $C$ is independent of $\sigma$ and $v_1$.
Let $M = 2 C |X_0|_{\Sigma}$  and choose  the  stopping time
\begin{align} \label{sigmav}
    \sigma_{v_1} := \inf \{t>0, 2^\a C^\a  |X_0|^{\a-1}_{\Sigma} \wt{\ve}^{\frac 1\theta}_t(v_1) >1 \} \wedge 1.
\end{align}
It follows that $\wt{\Phi} (\wt{\calx}_{\sigma_{v_1}}^M) \subset \wt{\calx}_{\sigma_{v_1}}^M$ and
$\wt{\Phi}$ is a contraction in $L^{\9}(0,\sigma_{v_1}; L^2) \cap L^q(0,\sigma_{v_1}; L^p)$.
Hence, similar arguments as in the proof of \cite[Theorem 2.1]{BRZ16}
yield that there exists a unique $H^1$-solution $\wt{z}=\wt{z}_{v_1}$ to \eqref{equa-wtz-damp} on $[0,\sigma_{v_1}]$.
In particular,  $\wt{z}_{v_1}$ exists on $[0,1]$ if $\sigma_{v_1}  = 1$.

Thus, in order to show that $\wt{z}_{v_1}$ exists on $[0,1]$ with high probability,
it suffices to prove that $\sigma_{v_1} = 1$ with high probability.

For this purpose, we consider
\begin{align} \label{wte}
   \wt{\ve}(v_1)
   :=& |h e^{(\a-1) \Re \wt{\vf}}|_{L^\theta(0,1)}^\theta
    = \int_0^\9 (1+s)^{- \frac{d(\a-1)-4}{2} \theta -2}  e^{(\a-1)\theta \Re \vf(s)} ds \nonumber \\
    =& \int_0^\9 (1+s)^{- \frac{d(\a-1)-4}{2} \theta -2} \prod \limits_{k=1}^N e^{(\a-1)\theta \Re \vf_k(s)} ds,
\end{align}
where
$$\Re \vf_k(t) := \int_0^t \Re v_k g_k(s) d\beta_k(s) -  \int_0^t (\Re v_k g_k(s))^2 ds.$$

Note that,
\begin{align} \label{C1}
   C_1 := \sup\limits_{t>0} \prod\limits_{k=2}^N \exp \{ (\a-1)\theta \Re \vf_k(t) \} <\9,\ \ a.s..
\end{align}
To this end, by the theorem on time change for continuous martingales (see e.g. \cite[Section 3.4]{KS91}),
there exists a Brownian motion $\wt{\beta}_k$ such that
$\bbp \circ (\Re v_k \int_0^\cdot g_k d\beta_k)^{-1} = \bbp \circ \wt{\beta}^{-1}_k((\Re v_k)^2 \int_0^\cdot g_k^2 ds)$.
Moreover, the law of the iterated logarithm for Brownian motion (see e.g. \cite[Section 2.9]{KS91}) implies that
$\lim_{t\to \9} \wt{\beta}_k(t)-t = -\9$, a.s..
Taking into account $\int_0^t g_k^2 ds \geq c_0^2 t \to \9$ as $t\to \9$, we have
\begin{align*}
   \bbp(\lim_{t\to \9} \Re \vf_k(t) = -\9)
   = \bbp(\lim_{t\to \9}  \wt{\beta}_k((\Re v_k)^2 \int_0^t g_k^2 ds) - (\Re v_k)^2 \int_0^t g_k^2 ds  = -\9) =1,
\end{align*}
which implies \eqref{C1}, as claimed.

For $k=1$, since $\inf_{t\geq 0} g_1(t)\geq c_0>0$, for any $\Re v_1 \geq 1$,
$(\Re v_1)^2 \int_0^t g_1^2(s) ds \geq c^2_0 t \to \9$ as $t\to \9$ a.s..
Then, by the law of the iterated logarithm,
there exist $c, t_*>0$, such that for any $t\geq t_*>0$ and any $\Re v_1 \geq 1$,
\begin{align*}
   \bigg|\Re v_1 \int_0^t g_1(s)d \beta_1(s) \bigg|
   \leq& c \((\Re v_1)^2 \int_0^t g^2_1(s)ds \ln\ln \((\Re v_1)^2 \int_0^t g^2_1(s)ds\)^{-1} \)^{\frac 12} \\
   \leq&   \frac 12 (\Re v_1)^2 \int_0^t g^2_1(s)ds,
\end{align*}
which implies that
\begin{align*}
   \Re \vf_1(t) \leq - \frac 12 (\Re v_1)^2 \int_0^t g_1^2(s) ds \leq -\frac 12 (\Re v_1)^2 c_0^2 t, \ \ t\geq t_*.
\end{align*}
Thus, using the dominated convergence theorem we have
\begin{align} \label{vf1-infi-0}
   &\int_{t_*}^\9 (1+s)^{- \frac{d(\a-1)-4}{2} \theta -2}  e^{(\a-1)\theta \Re \vf_1(s)} ds \nonumber \\
   \leq{}&  \int_{t_*}^\9 (1+s)^{- \frac{d(\a-1)-4}{2} \theta -2}  e^{-\frac 12 (\a-1)\theta c_0^2 (\Re v_1)^2s} ds \to 0,
\end{align}
as $\Re v_1 \to \9$, $\bbp$-a.s..

Moreover, $\mathbb{P}-a.s.$,
for any $0<s\leq t^*$,
we see that
$\wt{\beta}_1 (\int_0^s (\Re v_1)^2 g_1^2 dr) - \int_0^s (\Re v_1)^2 g_1^2 dr \to -\9$,
as $\Re v_1 \to \9$,
which implies that
\begin{align*}
   \wt{f}_{v_1}(s)
   := \exp\((\a-1)\theta \(\wt{\beta}_1 (\int_0^s (\Re v_1)^2 g_1^2 dr) - \int_0^s (\Re v_1)^2 g_1^2 dr \)\) \to 0, \ \ \Re v_1 \to \9.
\end{align*}
Taking into account
\begin{align*}
   \sup\limits_{\Re v_1 \geq 1} \sup\limits_{0<s\leq t^*} \wt{f}_{v_1}(s)
   \leq \sup\limits_{0<t<\9} \exp \{ (\a-1)\theta (\wt{\beta}_1(t)-t) \} \leq C <\9,
\end{align*}
we apply the dominated convergence theorem to obtain that
\begin{align*}
  \int_0^{t_*} (1+s)^{- \frac{d(\a-1)-4}{2} \theta -2} \wt{f}_{v_1}(s)  ds \to 0,\ \ as\ \Re v_1 \to \9,\ \ a.s.,
\end{align*}
which along with the theorem on time change for continuous martingales implies that
\begin{align} \label{vf1-t*-0}
   & \mathbb{P} \( \lim\limits_{\Re v_1 \to \9} \int_0^{t_*} (1+s)^{- \frac{d(\a-1)-4}{2} \theta -2} e^{(\a-1)\theta \Re \vf_1(s)} ds =0 \) \nonumber \\
  =& \mathbb{P} \( \lim\limits_{\Re v_1 \to \9} \int_0^{t_*} (1+s)^{- \frac{d(\a-1)-4}{2} \theta -2} \wt{f}_{v_1} (s) ds =0 \) =1.
\end{align}

Hence, it follows from \eqref{vf1-infi-0} and \eqref{vf1-t*-0} that, as $\Re v_1 \to \9$,
\begin{align} \label{vf1-0}
    \int_0^{\9} (1+s)^{- \frac{d(\a-1)-4}{2} \theta -2}  e^{(\a-1)\theta \Re \vf_1(s)} ds \to 0,\ \ a.s..
\end{align}

Thus, plugging \eqref{C1} and \eqref{vf1-0} into \eqref{wte} we get that a.s. for any $t\geq 0$,
\begin{align} \label{ve-v1}
  \wt{\ve}_t(v_1)
  \leq \wt{\ve}(v_1)
  \leq& C_1 \int_0^\9 (1+s)^{- \frac{d(\a-1)-4}{2} \theta -2} e^{(\a-1) \theta \Re \vf_1(s)} ds \to 0,\ \Re v_1 \to \9, \ a.s..
\end{align}
Then, in view of the definition of $\sigma_{v_1}$ in \eqref{sigmav},
we obtain
\begin{align} \label{sigmav-1}
   \bbp (\sigma_{v_1} = 1)  \to 1,\ \ as\ \Re v_1 \to \9,
\end{align}
which implies that $\wt{z}_{v_1}$ exists in the energy space $H^1$ on the interval $[0,1]$ with high probability
if $\Re v_1$ is sufficiently large.

Now, we consider all random variables being evaluated at $\omega \in \{\sigma_{v_1}= 1\}$.
By \eqref{equa-wtz-damp},  for all $t\in [0,1]$,
\begin{align*}
   \int |x|^2 |\wt{z}_{v_1}(t)|^2 dx
   = \int |x|^2 |X_0|^2 dx
      + 4 \int_0^t \Im \int x\cdot \na \ol{\wt{z}_{v_1}}(s) \wt{z}_{v_1}(s) dx ds.
\end{align*}
Then, by Cauchy's inequality,
\begin{align*}
    \int |x|^2 |\wt{z}_{v_1}(t)|^2 dx
   \leq& |X_0|^2_{\Sigma}
      + 4 \int_0^t \(\int |x|^2 |\wt{z}_{v_1}(s)|^2 dx \)^\frac 12 |\na \wt{z}_{v_1}|_2 ds \\
   \leq&   (|X_0|^2_{\Sigma} + 4 \|\wt{z}_{v_1}\|^2_{C([0,1];H^1)})
      + 4 \int_0^t \int |x|^2 |\wt{z}_{v_1}(s)|^2 dx  ds,
\end{align*}
which implies by Gronwall's inequality, because $\|\wt{z}_{v_1}\|_{C([0,1]; H^1)}<\9$, that
\begin{align*}
   \sup\limits_{t\in [0,1]} \int |x|^2 |\wt{z}_{v_1}(t)|^2 dx  <\9 .
\end{align*}
In particular, $\wt{z}_{v_1}(1) \in \Sigma$.

Therefore,
as in the proof of Theorem \ref{Thm-sca-X} $(ii)$,
by virtue of the equivalence between
the asymptotics of $\wt{z}_{v_1}$ at time $1$ and $z_{v_1}$ at infinity,
we conclude that
the original solution $z_{v_1}$ scatters at infinity in $\Sigma$,
i.e., $\lim_{t\to \9}|e^{it\Delta} z_{v_1}(t) -u_+|_{\Sigma} =0$ for some $u_+\in \Sigma$.
Thus,
\begin{align*}
   \{\sigma_{v_1} = 1\} \subset A_{v_1},
\end{align*}
where $A_{v_1}$ denotes the event that the solution $X$ to \eqref{equa-x} exists globally and scatters at infinity in $\Sigma$.
Taking into account \eqref{sigmav-1} we obtain \eqref{sca-damped-X} in the pseudo-conformal space. \\

The noise effect on scattering in the energy space can be proved similarly.
Taking into account \eqref{equa-z} and that $A(t)=-i\Delta$,
we set  $\calx_\tau^M := \{w\in L^\9(0,\tau; L^2) \cap L^q(0,\tau; L^p): \|w\|_{L^\9(0,\tau; H^1)} + \|w\|_{L^q(0,\tau; W^{1,p})} \leq M \}$
with $(p,q)$ as above
and define $\Phi$ on $\calx_\tau^M$  by
\begin{align*}
    \Phi(w) = e^{-it \Delta} X_0 - \lbb i \int_0^t e^{-i(t-s)\Delta} e^{-\vf(s)} F(e^{\vf(s)} w(s)) ds, \ \ w\in \calx_\tau^M.
\end{align*}
Similar  to \eqref{Phi-v.1} and \eqref{Phi-v.2},  for any $w_j \in \calx_\tau^M$, $j=1,2$,
if $\ve_\tau(v_1):= |e^{(\a-1) \Re \vf} |^\theta_{L^\theta(0,\tau)}$,
\begin{align*}
      \|\Phi(w_j)\|_{L^\9(0,\tau; H^1)} + \|\Phi(w_j)\|_{L^q(0,\tau; W^{1,p})}
      \leq   C|X_0|_{H^1} +  C \ve^{\frac 1 \theta}_\tau(v_1) M^\a,
\end{align*}
and
\begin{align*}
      &\|\Phi(w_1) - \Phi(w_2)\|_{L^\9(0,\tau; L^2)} +  \|\Phi(w_1) - \Phi(w_2)\|_{L^q(0,\tau; L^p)} \\
      \leq& C \ve_\tau^{\frac 1 \theta}(v_1)   M^{\a-1} \|w_1 -w_2\|_{L^q(0,\tau; L^p)},
\end{align*}
where $C$ is independent of $\tau$ and $v_1$.
Then, taking $M = 2 C |X_0|_{H^1}$  and
\begin{align} \label{tauv}
    \tau_{v_1} = \inf \{t>0, 2^\a C^\a  |X_0|^{\a-1}_{H^1} \ve^{\frac 1\theta}_t(v_1) >1 \},
\end{align}
and using similar arguments as in the previous case
we see that  there exists a unique $H^1$-solution $z (=z_{v_1})$ to \eqref{equa-z} on $[0,\tau_{v_1})$.
In particular, $z_{v_1}$ exists globally if $\tau_{v_1} = \9$.

Note that, similar  to \eqref{ve-v1}, for any $t>0$,
\begin{align} \label{tauv*}
  \ve_t(v_1)
  \leq \ve(v_1)
   := |e^{(\a-1) \Re \vf} |^{\theta}_{L^\theta(0,\9)}
   \to 0,
\end{align}
as $\Re v_1\to \9$, $\bbp$-a.s.,
which along with \eqref{tauv} implies that
\begin{align} \label{tv-1}
      \bbp (\tau_{v_1}=\9)
      \to 1,\ \ as\ \Re v_1\to \9.
\end{align}

Below we consider $\omega\in \{\tau_{v_1}=\9\}$.
As in \eqref{Phi-v.1},  for any $t\in (0,\9)$,
\begin{align} \label{z-vev}
  \|z_{v_1}\|_{L^q(0,t;W^{1.p})} \leq C |X_0|_{H^1} + C \ve^{\frac 1  \theta} (v_1)   \|z_{v_1}\|^\a_{L^q(0,t;W^{1,p})},
\end{align}
where $C$ is independent of $t$ and $v_1$.
In view of \eqref{tauv*}, choosing $\Re v_1$
large enough, if necessary, and
using \cite[Lemma 6.1]{BRZ16.2} we obtain
$\|z_{v_1}\|_{L^q(0,t;W^{1.p})} \leq C  < \9$,
with $C$ independent of $t$ and $v_1$.
Taking $t\to \9$ we get
\begin{align} \label{z-vev*}
  \|z_{v_1}\|_{L^q(0,\9;W^{1.p})}<\9.
\end{align}
Then,  similar  to \eqref{diff-e-u},
by virtue of \eqref{z-vev*}
we have that for $\Re v_1$ large enough,
\begin{align*}
   |e^{it_2\Delta} z_{v_1}(t_2) - e^{it_1\Delta}z_{v_1}(t_1)|_{H^1}
   \leq& C \|  e^{-\vf} F(e^{ \vf}z_{v_1}) \|_{L^{q'}(t_2,t_1;W^{1,p'})} \\
   \leq&  C({v_1}) \ve^{\frac {1}{\theta}}(v_1) \|z_{v_1}\|^{\a}_{L^{q}(t_2,t_1;W^{1,p})}
   \to  0,\ \ as\ t_1,t_2\to \9,
\end{align*}
which implies that  there exists $u_+ \in H^1$   such that
\begin{align*}
     e^{it\Delta} z_{v_1}(t) \to u_+,\ \ in\ H^1,\ as\ t \to \9,
\end{align*}
thereby yielding that $X$ scatters at  infinity in $H^1$.

Therefore,  for $\Re v_1 $ large enough,
\begin{align*}
    \{\tau_{v_1}=\9\}   \subset A_{v_1},
\end{align*}
which along with \eqref{tv-1}
implies \eqref{sca-damped-X}
in the energy space $H^1$.
This completes the proof. \hfill $\square$

\section{Strichartz and local smoothing estimates} \label{Sec-Stri}

In this section, we summarize the Strichartz and local smoothing estimates
used in this paper, mainly based on the work \cite{MMT08}.

Let $D_0 = \{|x| \leq 2\}$, $D_j = \{2^j\leq |x| \leq 2^{j+1}\}$,
and $D_{<j} = \{|x|\leq 2^j\}$, $j\geq 1$.
Set
$A_j = \bbr \times D_j$, $j\geq 0$,
and $A_{<j} = \bbr \times D_{<j}$, $j\geq 1$.
The local smoothing space is the completion of the Schwartz space
with respect to the norm
$\|u\|_{LS}^2 = \sum_{k=-\9}^\9 2^k \|S_k u\|^2_{{LS}_k}$,
and the dual norm is
$\|u\|_{LS'}^2 = \sum_{k=-\9}^\9 2^{-k} \|S_k u\|^2_{LS'_k}$,
where
$\{S_k\}$ is a dyadic partition of unity of frequency,
\begin{align*}
    &\|u\|_{LS_k} = \|u\|_{L^2_{t,x}(A_0)} + \sup\limits_{j>0} \|\<x\>^{-\frac 12} u\|_{L^2_{t,x}(A_j)},\ \ k\geq 0, \\
    &\|u\|_{LS_k} = 2^{\frac k2} \|u\|_{L^2_{t,x}(A_{<-k})} + \sup\limits_{j\geq -k} \|(|x|+2^{-k})^{-\frac12} u\|_{L^2_{t,x}(A_j)},\ \ k<0.
\end{align*}
(Since the notation $X$ stands for the solution to \eqref{equa-x},
in order to avoid confusions,
we use the different notation $LS$, instead of $X$ in \cite{MMT08}, for the local smoothing space.)

Similarly, for every $-\9\leq S < T \leq \9$,
we can also define $A_j$, $A_{<j}$ on the time interval $(S,T)$,
and $LS{(S,T)}$ denotes the local smoothing space defined on $(S,T)$.

We say that $(p,q)$ is a Strichartz pair,
if $2/q = d(1/2 - 1/p)$,
$(p,q) \in [2,\9] \times [2,\9]$,
and $(p,q,d) \not= (\9,2,2)$.

We first present the Strichartz and local smoothing estimates essentially proved in \cite{MMT08}.
\begin{theorem}   \label{Thm-Stri}
Consider the equation
\begin{align} \label{equa-stri}
   i\partial_t u = (\Delta + b\cdot \na + c) u + f
\end{align}
with $u(0)= u_0$,
$d\geq 3$. Assume that the coefficients $b,c$ satisfy
\begin{align} \label{cond-b}
   &\sum\limits_{j\in \mathbb{N}} \sup\limits_{A_j}
    \<x\> |b(t,x)| \leq \kappa,
\end{align}
and
\begin{align}
   &\sup\limits_{\bbr \times \bbr^d} \<x\>^2 (|c(t,x)|+|{\rm div}\ b(t,x)|) \leq \kappa,  \label{bc-ka}  \\
   &\limsup\limits_{|x|\to \9} \<x\>^2 (|c(t,x)| + |{\rm div}\ b(t,x)|) <\ve  \ll 1. \label{cond-bc}
\end{align}

$(i)$. For any $u_0\in L^2$, $T\in (0,\9)$ and any  two Strichartz pairs $(p_k,q_k)$, $k=1,2$,
we have the local-in-time Strichartz estimates, i.e.,
\begin{align}  \label{L-Stri}
    \|u\|_{L^{q_1}(0, T;L^{p_1}) \cap LS{(0,T)}}\leq
    C_T (|u_0|_{2}+\|f\|_{L^{q_2'}(0,T;L^{p_2'})+ LS'{(0,T)}}).
\end{align}

$(ii)$
Assume in addition that \eqref{cond-b} and \eqref{bc-ka} hold also
for $\p_j b$ and $\p_j c$, $1\leq j\leq d$.
Then, for any Strichatz pairs $(p_k, q_k)$, $1\leq k\leq 3$,
\begin{align}  \label{L-Stri*}
    \|\partial_ju\|_{L^{q_1}(0,T;L^{p_1}) \cap LS{(0,T)}}
    \leq&
    C_T \big(|u_0|_{H^1}+\|f\|_{L^{q_2'}(0,T;L^{p_2'})+ LS'{(0,T)}} \nonumber \\
       &\qquad+ \|\partial_j f\|_{L^{q_3'}(0,T;L^{p_3'})+ LS'{(0,T)}} \big),
\end{align}

$(iii)$. Assume in addition that $\kappa \leq \ve \ll 1$. Then,
we have the global-in-time Strichartz estimates, i.e.,
for any Strichartz pairs $(p_k,q_k)$, $k=1,2$,
\begin{align}  \label{G-Stri}
    \|u\|_{L^{q_1}(\bbr;L^{p_1}) \cap LS}\leq
    C (|u_0|_{2}+\|f\|_{L^{q_2'}(\bbr;L^{p_2'})+ LS'}).
\end{align}

\end{theorem}

{\bf Proof.}
$(i)$. The proof is similar to that of \cite[Lemma 4.1]{BRZ14}.

$(ii)$.  Estimate \eqref{L-Stri*} can be proved similarly as in the proof of
\cite[Lemma 2.7]{BRZ16}.
In fact, for each $1\leq j\leq d$,   $v_j:= \partial_j u$ satisfies
\begin{align} \label{equa-vj}
    i \partial_t v_j
    = (\Delta + b\cdot \na + c) v_j
       + (\partial_j b \cdot \na + \partial_j c) u
       + \partial_j f
\end{align}
with $v_j(0) = \partial_j u_0$.
Then, applying \eqref{L-Stri} to \eqref{equa-vj}
and using
\begin{align} \label{X'-X}
      \|(b\cdot \na + c) u \|_{LS'} \leq   C \kappa \|u\|_{LS}
\end{align}
(see \cite[Proposition 2.3]{MMT08};
note that the space $\wt{X}$ in \cite[Proposition 2.3]{MMT08} coincides with the local smoothing space $LS$ defined above if $d\geq 3$)
we obtain
\begin{align*}
   & \|v_j\|_{L^{q_1}(0,T;L^{p_1}) \cap LS{(0,T)}}\\
    \leq{}&
       C_T (|\partial_j u_0|_{2}+ \|(\partial_j b\cdot \na + \partial_j c) u\|_{LS'{(0,T)}} +\|\partial_j f\|_{L^{q_3'}(0,T;L^{p_3'})+ LS'{(0,T)}}) \\
    \leq{}&  C_T (|u_0|_{H^1}+ \|u\|_{LS{(0,T)}} +\|\partial_j f\|_{L^{q_3'}(0,T;L^{p_3'})+ LS'{(0,T)}}).
\end{align*}
Thus, applying again \eqref{L-Stri} yields immediately \eqref{L-Stri*}.

$(iii)$. Since for the Laplacian $-\Delta$  in dimension $d\geq 3$
the associated bicharacteristic flow is not trapped
and zero is not an eigenvalue or a resonance,
we use \cite[Theorem 1.22]{MMT08} and \eqref{X'-X} to obtain that
\begin{align*}
      \|u\|_{L^{q_1}(\bbr;L^{p_1}) \cap LS}
    \leq& C (|u_0|_{2}+ \|(b\cdot \na + c) u\|_{LS'} +\|f\|_{L^{q_2'}(\bbr;L^{p_2'})+ LS'}) \\
    \leq& C (|u_0|_{2}+  \ve \|u\|_{LS} +\|f\|_{L^{q_2'}(\bbr;L^{p_2'})+ LS'})
\end{align*}
Thus, taking $\ve$ small enough we prove \eqref{G-Stri}. \hfill $\square$

\begin{remark} \label{Rem-Stri}
Using characteristic functions we see that,
the estimates \eqref{L-Stri}-\eqref{G-Stri} and \eqref{X'-X} are also valid on $(S,T)$ for any $-\9 \leq S<T \leq \9$
if the corresponding conditions hold on $(S,T)$.
\end{remark}

\begin{corollary} \label{Cor-Stri}
$(i)$.  Consider  \eqref{equa-stri} and assume that for any multi-index $0\leq |\beta|\leq 1$,
\begin{align}
    & \sup\limits_{(t,x)\in \bbr^+ \times \bbr^d}
   \<x\>^2 (|\partial_x^\beta b(t,x)| + |c(t,x)|)  <\9,  \label{cond-bc-0}
\end{align}
and
\begin{align}
    &\limsup \limits_{t\to \9} \sup\limits_{x\in \bbr^d} \<x\>^2 (|\partial_x^\beta b(t,x)| + |c(t,x)|) \leq \ve \ll 1, \label{cond-bc-2}  \\
    &  \limsup\limits_{|x|\to \9} \sup\limits_{0\leq t<\9} \<x\>^2 (|\partial_x^\beta b(t,x)| + |c(t,x)|)  \leq \ve \ll 1.   \label{cond-bc-1}
\end{align}
Then, for any $ 0\leq S < T \leq \9$ and any two Strichartz pairs $(p_k,q_k)$, $k=1,2$,
we have the global-in-time Strichartz estimates
\begin{align}  \label{Cor-Stri-esti}
    \|u\|_{L^{q_1}(S,T;L^{p_1}) \cap LS{(S,T)}}\leq
    C (|u(S)|_{2}+\|f\|_{L^{q_2'}(S,T;L^{p_2'})+ LS'{(S,T)}}),
\end{align}
where $C$ is independent of $S,T$.

$(ii)$. Assume in addition that, for every $1\leq j\leq d$,
$\partial_j b$ and $\partial_j c$ also satisfy \eqref{cond-bc-2} and \eqref{cond-bc-1}.
Then, for any Strichartz pairs $(p_k,q_k)$, $1\leq k\leq 3$,
\begin{align}  \label{Cor-Stri-esti-H1}
    \|\partial_ju\|_{L^{q_1}(S,T;L^{p_1}) \cap LS{(S,T)}}
    \leq& C (|u(S)|_{H^1}+\|f\|_{L^{q_2'}(S,T;L^{p_2'})+ LS'{(S,T)}} \nonumber \\
        &\qquad + \|\partial_j f\|_{L^{q_3'}(S,T;L^{p_3'})+ LS'{(S,T)}}),
\end{align}
where $C$ is independent of $S,T$.

$(iii)$. Assume the conditions in  $(i)$ (resp. $(ii)$) above to hold.
Then, \eqref{Cor-Stri-esti} (resp. \eqref{Cor-Stri-esti-H1}) also holds
with $u_0$ replaced by the final datum $u(T)$.
\end{corollary}

{\bf Proof.} $(i)$. It suffices to prove the assertion for $S=0$ and $T=\9$.
By \eqref{cond-bc-1}, for any $T_1>0$,
\begin{align}
   &  \sup\limits_{(t,x)\in [0,T_1] \times \bbr^d}
   \<x\>^2 (|\partial_x^\beta b(t,x)| + |c(t,x)|) \leq \kappa(T_1) <\9, \label{bcT.1} \\
   &\limsup\limits_{|x|\to \9}\sup\limits_{t\in [0,T_1]}
   \<x\>^2 (|\partial_x^\beta b(t,x)| + |c(t,x)|) \leq \ve \ll 1,  \label{bcT.1*}
\end{align}
which implies that \eqref{cond-b}-\eqref{cond-bc} hold in
the time range $[0,T_1]$.
Thus, using Theorem \ref{Thm-Stri} $(i)$ (see also Remark \ref{Rem-Stri})
we obtain
\begin{align} \label{bcT.1**}
    \|u\|_{L^{q_1}(0,T_1;L^{p_1}) \cap LS{(0,T_1)}}\leq
    C_{T_1} (|u_0|_{2}+\|f\|_{L^{q_2'}(0,T_1;L^{p_2'})+ LS'{(0,T_1)}}).
\end{align}

Moreover, for $T_1$ fixed and large enough, by \eqref{cond-bc-2} and \eqref{cond-bc-1},
\begin{align}
   &\sup\limits_{(t,x)\in [T_1, \9) \times \bbr^d}
   \<x\>^2 (|\partial_x^\beta b(t,x)| + |c(t,x)|) \leq 2 \ve  \ll 1,  \label{bcT.2} \\
   &\limsup\limits_{|x|\to \9}\sup\limits_{t\in [T_1,\9)}
   \<x\>^2 (|\partial_x^\beta b(t,x)| + |c(t,x)|) \leq 2 \ve  \ll 1, \label{bcT.2*}
\end{align}
which implies that \eqref{cond-b}-\eqref{cond-bc} hold on $[T_1,\9)$.
Thus, using Theorem \ref{Thm-Stri} $(iii)$ we get
\begin{align} \label{bcT.2**}
    \|u\|_{L^{q_1}(T_1,\9;L^{p_1}) \cap LS{(T_1,\9)}}
    \leq& C (|u(T_1)|_{2}+\|f\|_{L^{q_2'}(T_1,\9;L^{p_2'})+ LS'{(T_1,\9)}}) \nonumber \\
    \leq& C_{T_1}  (|u_0|_{2}+\|f\|_{L^{q_2'}(0,\9;L^{p_2'})+ LS'{(0,\9)}}),
\end{align}
where in the last step we also used \eqref{bcT.1**} with the Strichartz pair $(p_1,q_1)=(2,\9)$,
and $C$ is independent of $T_1$.
Combining \eqref{bcT.1**} and \eqref{bcT.2**} we obtain \eqref{Cor-Stri-esti}.

$(ii)$. The argument is similar to that in the proof of Theorem \ref{Thm-Stri} $(ii)$.

$(iii)$. It is sufficient to consider $S=0$.
Let $T_1$ be as above and assume $T>T_1$ without loss of generality.
Set $\wt{g}(t)= g(T-t)$, where $g= u, b, c, f$, $t\in[0,T]$.
Then, by \eqref{equa-stri},
\begin{align*}
   i \p_t \wt{u} = (-1) (\Delta + \wt{b}(t)\cdot \na + \wt{c}(t)) \wt{u} - \wt{f}.
\end{align*}

Note that,
\eqref{bcT.2} and \eqref{bcT.2*} hold for $\wt{b}$ and $\wt{c}$ with $[T_1,\9)$ replaced by $[0,T-T_1]$.
Applying Theorem \ref{Thm-Stri} $(iii)$ on $[0,T-T_1]$,
we get
\begin{align}\label{wtu-T1.1}
         \|\wt{u}\|_{L^{q_1}(0,T-T_1;L^{p_1}) \cap LS{(0,T-T_1)}}
    \leq& C (|\wt{u}(0)|_{2} +\|\wt{f}\|_{L^{q_2'}(0,T-T_1;L^{p_2'})+ LS'{(0,T-T_1)}}) \nonumber \\
    \leq& C (|u(T)|_{2}+ \|f\|_{L^{q_2'}(T_1,T;L^{p_2'})+ LS'{(T_1,T)}}),
\end{align}
where $C$ is independent of $T$.

Moreover, since \eqref{bcT.1} and \eqref{bcT.1*} hold for $\wt{b}, \wt{c}$ on $[T-T_1, T]$ replacing $[0,T_1]$,
applying Theorem \ref{Thm-Stri} $(i)$ on $[T-T_1,T]$ we have
\begin{align} \label{wtu-T1.2}
         \|\wt{u}\|_{L^{q_1}(T-T_1,T;L^{p_1}) \cap LS{(T-T_1,T)}}
    \leq& C_{T_1} (|\wt{u}(T-T_1)|_{2} +\|\wt{f}\|_{L^{q_2'}(T-T_1,T;L^{p_2'})+ LS'{(T-T_1,T)}}) \nonumber \\
    \leq&  C_{T_1} (|u(T)|_{2}  +\|f\|_{L^{q_2'}(0,T;L^{p_2'})+ LS'{(0,T)}})
\end{align}
with  $C_{T_1}$  independent of $T$,
where in the last step we used \eqref{wtu-T1.1} with $(p_1,q_1)=(2,\9)$.

Therefore, putting together \eqref{wtu-T1.1} and \eqref{wtu-T1.2} we prove \eqref{Cor-Stri-esti}
with $u_0$ replaced by $u(T)$.
The proof for $\p_j u$ is similar. \hfill $\square$ \\

In the remainder of this section we verify the global-in-time
Strichartz and local smoothing estimates used in Sections \ref{Sec-Finite-Varia}--\ref{Sec-Sca-Damp}.

First, consider the
global-in-time Strichartz and local smoothing estimates for the operator $\wt{A}_*$ on $[0,1)$
in Section \ref{Sec-Finite-Varia}.
We take $\partial_{jh} \wt{b}_* (t,\xi)$ for an example
to verify the conditions \eqref{cond-bc-0}-\eqref{cond-bc-1}
under Assumptions $(H0)$ and $(H1)$, $1 \leq j,h\leq d$.

By \eqref{b*} and \eqref{wtb},
\begin{align*}
   \partial_{jh} \wt{b}_*(t,x)
   =&- 2 (1-t)^{-3} \sum\limits_{k=1}^N \na \partial_{jh} \phi_k(\frac{x}{1-t}) \int_{\frac{t}{1-t}}^\9 g_k(s)d\beta_k(s) \\
     &+2 (1-t)^{-3} \sum\limits_{k=1}^N \na \partial_{jh}((\Re\phi_k)\phi_k)(\frac{x}{1-t}) \int_{\frac{t}{1-t}}^\9 g_k^2 (s)ds.
\end{align*}
Then, in view of the asymptotic flatness condition \eqref{AF-3},
we obtain that
\begin{align} \label{jhb.1}
  \<x\>^2 |\partial_{jh} \wt{b}_*(t,x)|
  \leq& 2 \sum\limits_{k=1}^N \<x\>^2 \( \big|\na \partial_{jh} \phi_k(\frac{x}{1-t})| + |\na \partial_{jh}((\Re\phi_k)\phi_k)(\frac{x}{1-t}) \big| \) \nonumber \\
      & \qquad \cdot (1-t)^{-3} \( \bigg|\int_{\frac{t}{1-t}}^\9 g_k(s)d\beta_k(s)\bigg|+ \int_{\frac{t}{1-t}}^\9 g_k^2 (s)ds \) \nonumber \\
  \leq& 2\sum\limits_{k=1}^N \wt{\ve}_k(\frac{x}{1-t}) \wt{r}_k(t) ,
\end{align}
where $\wt{\ve}_k(x) \to 0$ as $|x|\to \9$, and $\wt{r}_k(t) := (1-t)^{-3} ( |\int_{\frac{t}{1-t}}^\9 g_k(s)d\beta_k(s)|+ \int_{\frac{t}{1-t}}^\9 g_k^2 (s)ds )$.
Note that
\begin{align} \label{jhb.2}
     \sup\limits_{x\in \bbr^d} \sup\limits_{0\leq t<1} \wt{\ve}_k(\frac{x}{1-t}) <\9, \ \
     \lim\limits_{|x|\to \9} \sup\limits_{0\leq t<1} \wt{\ve}_k(\frac{x}{1-t}) =0.
\end{align}
Moreover, by the law of the iterated logarithm,
\begin{align*}
     \bigg|\int_{\frac{t}{1-t}}^\9 g_k(s)d\beta_k(s)\bigg|
     \leq  \( 2\int_{\frac{t}{1-t}}^\9 g_k^2 (s)ds \ln \ln \(\int_{\frac{t}{1-t}}^\9 g_k^2 (s)ds \)^{-1} \)^{\frac 12},\ \ a.s.,
\end{align*}
which along with \eqref{ILog} implies that
\begin{align} \label{jhb.3}
   \sup\limits_{0\leq t<1} \wt{r}_k(t) <\9,\ \
   \limsup\limits_{t\nearrow 1} \wt{r}_k(t) =0.
\end{align}
Thus, putting together \eqref{jhb.1}-\eqref{jhb.3} yields that the conditions
\eqref{cond-bc-0}-\eqref{cond-bc-1} hold for $\p_{jh} \wt{b}_*$
with $\9$ and $\ve$ replaced by $1$ and $0$, respectively.

Similar arguments also apply to $\partial_x^\g \wt{b}_*$ and $\partial_x^\g \wt{c}_*$, $0\leq |\g| \leq 1$,
which, via Corollary \ref{Cor-Stri} (see also Remark \ref{Rem-Stri}),
imply global-in-time
Strichartz and local smoothing estimates for $\wt{A}_*$ on $[0,1)$.  \\

Next, we check that global-in-time Strichartz and local smoothing estimates for the operator $A_*$ in \eqref{equa-z*} and \eqref{equa-v}
under the condition \eqref{asymflat} and that $g_k\in L^2(\bbr^+)$, a.s., $1\leq k\leq N$.

We illustrate this for $b_*$ only.
For any $0\leq |\beta| \leq 2$, by \eqref{b*},
\begin{align*}
   \partial_x^\beta b_*(t,x)
   = &(-2) \sum\limits_{k=1}^N \na \partial_x^\beta \phi_k(x) \int_t^\9 g_k(s) d\beta_k(s)\\
     &{}+ 2 \sum\limits_{k=1}^N  \na \partial_x^\beta ((\Re \phi_k) \phi_k)(x) \int_t^\9 g^2_k(s) ds,
\end{align*}
which implies that
\begin{align} \label{betab.1}
   |\<x\>^2 \partial_x^\beta b_*(t,x)|
   \leq& C  \sum\limits_{k=1}^N \<x\>^2(| \na \partial_x^\beta \phi_k(x)|+ |\na \partial_x^\beta  ((\Re \phi_k) \phi_k)(x)|)  \nonumber \\
       &  \qquad \cdot \(\bigg|\int_t^\9 g_k(s) d\beta_k(s) \bigg| + \int_t^\9 g^2_k(s) ds \) \nonumber \\
      =:& C  \sum\limits_{k=1}^N  \ve_k(x) r_k(t).
\end{align}
Note that, \eqref{asymflat} implies that
\begin{align} \label{betab.2}
    \sup\limits_{x\in \bbr^d} \ve_k(x)<\9,\ \
    \lim\limits_{|x|\to \9} \ve_k(x)=0.
\end{align}
Moreover, since $g_k\in L^2(\bbr^+)$, a.s.,
we have that $|\int_0^\9 g_k(s) d\beta_k(s)| <\9$, a.s., and
by the law of iterated logarithm,
as $t\to \9$,
\begin{align*}
   \bigg|\int_t^\9 g_k(s) d\beta_k(s) \bigg|
   \leq \(2 \int_t^\9 g_k^2 ds \ln \ln (\int_t^\9 g_k^2 ds) \)^{\frac 12}  \to 0,\ \ a.s..
\end{align*}
Hence,
\begin{align} \label{betab.3}
     \sup\limits_{0\leq t<\9} r_k(t) <\9,\ \
   \limsup\limits_{t\to \9} r_k(t) =0.
\end{align}
Thus, we conclude from \eqref{betab.1}-\eqref{betab.3} that
$\p_x^\beta b_*$,  $0\leq |\beta|\leq 2$,  satisfy the conditions \eqref{cond-bc-0}-\eqref{cond-bc-1} with $\ve=0$.
Similar arguments also apply to $\partial_x^\g c_*$, $0\leq |\g|\leq 1$.
Thus, global-in-time Strichartz and local smoothing estimates for $A_*$ follow from Corollary \ref{Cor-Stri}. \\

\section{Appendix}  \label{Sec-App}

\subsection{Proof of Theorem \ref{Thm-GWP}}

The proof is similar to that in \cite[Theorem 2.2]{BRZ14} and \cite[Theorem 1.2]{BRZ16}.

First note that, under Assumptions $(H0)$,
local-in-time Strichartz and local smoothing estimates hold for the operator $A(t)$ in \eqref{equa-z}.
In fact, by Assumption $(H0)$ and the Burkholder-Davis-Gundy inequality,
for any $0<T<\9$,
\begin{align*}
   \bbe \sup\limits_{t\in [0,T]}
   \bigg|\int_0^t g_k d\beta_k(s) \bigg|
   \leq C \bbe \(\int_0^T g^2_k ds\)^{\frac 12}
   \leq C \|g_k\|_{L^\9(\Omega \times (0,T))} T^\frac 12 <\9,
\end{align*}
which implies that
\begin{align*}
   C_T:=  \sup\limits_{t\in[0,T]} \bigg|\int_0^t g_k(s) d\beta_k(s) \bigg| + \int_0^T g_k^2(s)ds<\9,\ \ \bbp-a.s..
\end{align*}
Moreover, by \eqref{y-b}, for $0\leq |\beta|\leq 2$,
\begin{align*}
   |\<x\>^2 \partial_x^\beta b(t,x)|
   \leq& C  \sum\limits_{k=1}^N \<x\>^2(| \partial_x^\beta \na \phi_k(x)|+ |\partial_x^\beta \na ((\Re \phi_k) \phi_k)(x)|)  \\
       &  \qquad \cdot \(\bigg|\int_0^t g_k(s) d\beta_k(s)\bigg| + \int_0^t g^2_k(s) ds \).
\end{align*}
Thus, by \eqref{asymflat},
we see that $\p_x^\beta b$, $0\leq |\beta|\leq 2$, satisfy the conditions \eqref{cond-b}--\eqref{cond-bc}
on $[0,T]$ with $b$ and ${\rm div}\ b$ replaced by $\p_x^\beta b$.
The argument for the coefficients $\p_x^\g c$ is similar, $0\leq |\g|\leq 1$.
Therefore, by Theorem \ref{Thm-Stri} $(ii)$,  local-in-time Strichartz
and local smoothing estimates hold for $A(t)$ in \eqref{equa-z}.

Now, the local well-posedness for \eqref{equa-x} follows from similar arguments
as in the proof of  \cite[Proposition $2.5$]{BRZ16}.
As regards the global well-posedness, applying It\^o's formula we have $\bbp$-a.s.
for any $t\in [0,\tau^*)$,
where $\tau^*$ is the maximal existing time,  that
\begin{align} \label{mass}
  |X(t)|_2^2 = |X_0|_2^2 + \sum\limits_{k=1}^N \int_0^t \sigma_{0,k}(s) d\beta_k(s) ,
\end{align}
and for the Hamiltonian  $H(X) :=\frac 12 |\na X|_2^2 - \frac{\lbb}{\a+1} |X|^{\a+1}_{L^{\a+1}}$,
\begin{align} \label{hamil}
   H(X(t)) =H(X_0) + \int_0^t a_1(s) ds +  \sum\limits_{k=1}^N \int_0^t \sigma_{1,k}(s) d\beta_k(s),
\end{align}
where
\begin{align*}
   a_1(s)=& - \Re \int \na \ol{X}(s) \na(\mu(s) X(s)) dx
              + \frac 12 \int |\na (G_k(s) X(s))|^2 dx  \\
          &    - \frac{\lbb(\a-1)}{2} \int (\Re G_k(s))^2 |X(s)|^{\a+1} dx, \\
    \sigma_{0,k}(s) =& 2 \int \Re G_k(s) |X(s)|^2 dx, \\
    \sigma_{1,k}(s) =&  \Re \int \na \ol{X}(s) \na (G_k (s) X(s)) dx
              -\lbb  \int \Re G_k(s) |X(s)|^{\a+1} dx .
\end{align*}

Thus, similar arguments   as in the proof of \cite[Lemma $3.1$]{BRZ14} and \cite[Theorem $3.1$]{BRZ16}
imply that
$\sup_{[0,\tau^*)}|X(t)|_{H^1}^2 < \9$, a.s.,
thereby yielding that
$X$ exists globally.

Moreover, if $X_0\in \Sigma$, we have that for
the virial function $V(X(t)):= \int |x|^2 |X(t,x)|^2 dx$,
\begin{align} \label{Ito-V}
    V(X(t)) = V(X_0) + 4 \int_0^t G(X(s)) ds + \sum\limits_{k=1}^N \int_0^t \sigma_{2,k}(s) d\beta_k(s),
\end{align}
where
$G(X)= \Im \int x \cdot \na \ol{X} X d x$
and
$\sigma_{2,k}(s) := 2 \int \Re G_k(s) |x|^2 |X(s)|^2 dy$.
Then, arguing as in the proof of \cite[$(4.6)$]{BRZ14.3} we obtain $X\in L^2 (\Omega; C([0,T]; \Sigma))$.

The proof for the estimate \eqref{thmx2*} in weighted spaces is similar to that
of Lemma \ref{Lem-wtvT}.
Actually,
consider the Strichartz pair $(p,q)$ and $1<\theta<\9$ as in the proof of Lemma \ref{Lem-wtvT}.
We only need to prove that, for $z=e^{-\vf}X$ and for each $1\leq j\leq d$,
\begin{align} \label{xjz}
\|x_jz\|_{L^q(0,T; L^p)}<\9\ \ a.s..
\end{align}
Since $\mathbb{P}$-a.s. $X\in L^{q}(0,T; W^{1,p})$, $\vf \in C([0,T]; W^{1,\9})$,
and $b\in C([0,T]; L^\9)$,
we  take a finite partition $\{(t_k,t_{k+1})\}$ of $[0,T]$,
such that $\|z\|_{L^q(t_k,t_{k+1}; W^{1,p})} \leq \ve$.
Then, similar to \eqref{lpq-xwtz},
\begin{align*}
  \|x_jz\|_{L^q(t_k,t_{k+1}; L^p)}
  \leq C_T \|z\|_{C([0,T]; \Sigma)}
       + C_T \ve^{\a-1} \|x_jz\|_{L^q(t_k,t_{k+1}; L^p)},
\end{align*}
which implies \eqref{xjz} by taking $\ve$ sufficiently small and
summing over $k$.

Regarding \eqref{glo-X-H1},
we set $\calh(X):= \frac 12 |X|_{H^1}^2 - \frac{\lbb}{\a+1} |X|_{L^{\a+1}}^{\a+1}
= H(X) + \frac 12|X|_2^2$.
Applying It\^o's formula and using \eqref{mass} and \eqref{hamil}
we get that for any $m\geq 2$,
\begin{align*}
   (\calh(X(t)))^m
   =& (\calh(X(0)))^m
       + m \int_0^t (\calh(X(r)))^{m-1} a_1(X(r)) dr \\
    & + \frac 12 m(m-1) \sum\limits_{k=1}^N  \int_0^t (\calh(X(r)))^{m-2}
       (\frac 12 \sigma_{0,k}(r) + \sigma_{1,k}(r))^2 dr
      + \calm (t) \\
   \leq& C(m) |X_0|_{H^1}^m
        + C(m) \sum\limits_{k=1}^N  \int_0^t (\calh(X(r)))^m g_k^2(r)dr + \calm(t),
\end{align*}
where  $\calm(t) =\sum_{k=1}^N \int_0^t m (\calh(X(r)))^{m-1} (\frac 12 \sigma_{0,k}(r)+ \sigma_{1,k}(r)) d\beta_k(r)$.

Then, apply the stochastic Gronwall inequality in \cite[Lemma 2.2]{RZ18} we obtain that
for any $0<q<p<1$,
\begin{align*}
   \bbe \sup\limits_{0\leq s\leq t}
   (\calh(X(s)))^{mq}
   \leq C(m)
    |X_0|_{H^1}^{mq}
    (\frac{p}{p-q})
    \(\bbe e^{C(m) \frac{p}{1-p} \sum\limits_{k=1}^N \int_0^t g_k^2(s) ds}\)^{\frac{1-p}{p}q}.
\end{align*}
Then, taking $q= p/2\in (0,1)$,  using Fatou's lemma
we obtain that for any $m\geq 2$ and $p\in(0,1)$,
\begin{align} \label{esti-H-p}
   \bbe \sup\limits_{0\leq t< \9}
   (\calh(X(t)))^{\frac 12 mp}
   \leq 2C(m) |X_0|_{H^1}^{\frac 12 mp}
    \(\bbe e^{C(m)\frac{p}{1-p} \sum\limits_{k=1}^N \int_0^\9 g_k^2(t) dt}\)^{\frac 12 },
\end{align}
which implies \eqref{glo-X-H1}
because $m$ and $q$ can be chosen arbitrarily.

Therefore, the proof is complete.
\hfill $\square$

\subsection{Proof of Proposition \ref{Prop-wtu}}\label{sec:proof-pro}
The proof is similar to that of Corollary \ref{Cor-wtz-globH1}
based on the pseudo-conformal energy
and the global bound \eqref{glbdd-wtz-LpWq-Weight} of $\wt{z}_*$.

In fact,
in the case where $\a\in [1+4/d, 1+4/(d-2)]$,
similarly to \eqref{Ito-wtE1},
\begin{align*}
   \wt{E}_1(\wt{u}(t))
   =& \wt{E}_1(\wt{u}(\wt{T}))
     + \frac{16}{\a+1} (1- \frac{d}{4}(\a-1))
              \int_{\wt{T}}^t (1-r)^{\frac d2 (\a-1) -3} |\wt{X}(r)|_{L^{\a+1}}^{\a+1} dr \\
   \leq&  \wt{E}_1(z_*(\wt{T}))
   \leq  C (1+\|\wt{z}_*\|^{\a+1}_{C([0,1);H^1)})
   <\9,  \ \ t\in[\wt{T}, 1),\ a.s.,
\end{align*}
where the last step is due to
the uniform bound \eqref{glbdd-wtz-LpWq-Weight}
and $C$ is independent of $\wt{T}$.
This yields that for some positive constant $C$ independent of $\wt{T}$,
\begin{align*}
    \|\wt{u}\|^2_{C([\wt{T}, 1); H^1)} \leq C(1+\|\wt{z}_*\|^{\a+1}_{C([0,1);H^1)}),\ \  a.s..
\end{align*}

Moreover, in the case where $\a\in (1+\a(d), 1+4/d)$,
similarly to \eqref{Ito-wtE2},
for any $t\in[\wt{T}, 1)$,
\begin{align*}
  \wt{E}_2(\wt{u}(t))
  =& \wt{E}_2(\wt{u}(\wt{T}))
     -8  (1- \frac{d}{4}(\a-1)) \int_{\wt{T}}^t (1-r)^{1 - \frac d2 (\a-1)} |\na \wt{X}(r)|_2^2 dr \\
  \leq&  \wt{E}_2(\wt{z}_*(\wt{T}))
  \leq C (1+\|\wt{z}_*\|^{\a+1}_{C([0,1);H^1)})
  <\9, \  a.s.,
\end{align*}
which implies that for some positive constant $C$ independent of $\wt{T}$,
\begin{align*}
  \sup_{\wt{T}\leq t<1} |\wt{u}(r)|_{L^{\a+1}}^{\a+1}
    \leq C (1+\|\wt{z}_*\|^{\a+1}_{C([0,1);H^1)}),\ \  a.s..
\end{align*}
Thus,
using the arguments as those below \eqref{bdd-wtz-H1-case1} and \eqref{bdd-wtz-Lp-case2}
we obtain  \eqref{esti-wtz-LpWq-Weight-case1}
with $\wt{u}$ replacing $\wt{z}_*$.
Precisely,
for $\a \in (1+\a(d), 1+ 4/(d-2)]$
and any Strichartz pair $(\rho, \g)$,
\begin{align*}
\|\wt{u}\|_{L^{\gamma}(\wt{T},1;W^{1,\rho})}
     +\| |\cdot| \wt{u}\|_{L^{\gamma}(\wt{T},1;L^{\rho})}
    \leq C (|\wt{u}(\wt{T})|_\Sigma + \|\wt{u}\|_{C([0,1);H^1)})
    \leq C   <\9,
\end{align*}
where $C$ depends on $\|\wt{z}_*\|_{C([0,1);\Sigma)}$
and is independent of $\wt{T}$.
This also yields uniform bound for the $LS(\wt{T},1)$-norms
by Strichartz estimates.

Finally, using the uniform bound \eqref{glbdd-wtu-LpWq-Weight}
and the  local well-posedness arguments as in \cite{BRZ16}
we can extend the solution $\wt{u}$ to time $1$.
In particular, $\|\wt{u}\|_{C([0,1];\Sigma)} <\9$, a.s..

Therefore, the proof of Proposition \ref{Prop-wtu} is complete. \hfill $\square$

\subsection*{Acknowledgements}
D. Zhang  is supported by NSFC (No. 11501362).
Financial support by the DFG through CRC 1283 is gratefully acknowledged.

\end{document}